\documentclass[10pt]{IEEEtran}
\usepackage{graphicx,epstopdf,multirow,multicol,subfigure,booktabs}

\usepackage{amssymb,amsthm}
\usepackage{algorithm}
\usepackage[noend]{algpseudocode}
\usepackage{color}
\usepackage{pgfplotstable}
\usepackage{pgfplots}
\pgfplotsset{compat=newest}
\usepackage{float}
\usepackage{setspace}
\usepackage[cmex10]{amsmath}
\usepackage{bm,mathtools}
\usepackage{cite}
\usepackage{times}
\usepackage{upgreek}
\usepackage{enumitem}

\algnewcommand\algorithmicinput{\textbf{Input:}}
\algnewcommand\Input{\item[\algorithmicinput]}

\algnewcommand\algorithmicoutput{\textbf{Output:}}
\algnewcommand\Output{\item[\algorithmicoutput]}

\usepackage{hyperref}
\usepackage{booktabs}
\allowdisplaybreaks
\graphicspath{figures/}

\theoremstyle{remark}

\usepackage{float}
\floatstyle{ruled}
\newfloat{model}{thp}{lop}
\floatname{model}{Model}

\usepackage{amssymb,bm}

\newcommand{\omitit}[1]{}

\newcommand{\matpower}{\mbox{\sc Matpower}}
\newcommand{\pglib}{{\sc PGLib}}
\newcommand{\pglibopf}{{\sc PGLib-OPF}}

\usepackage{multirow}
\usepackage{varwidth}
\newcolumntype{M}[1]{>{\begin{varwidth}[t]{#1}}l<{\end{varwidth}}}

\usepackage[american,fulldiode]{circuitikz} %%Circuit drawing tools

\begin{document}

\title{The Power Grid Library for Benchmarking \\ AC Optimal Power Flow Algorithms}
\author{Contributors to the IEEE PES PGLib-OPF Task Force \\
\vspace{0.2cm}
{
\vspace{0.5cm}
\hspace{0.5cm}
\centering
\small
\begin{tabular}[t]{c@{\extracolsep{1cm}}c} 
Sogol Babaeinejadsarookolaee & Adam Birchfield \\
Electrical and Computer Engineering & Electrical and Computer Engineering \\ 
University of Wisconsin-Madison, WI, USA & Texas A\&M University, TX, USA \\
& \\
Richard D. Christie & Carleton Coffrin$^{*}$ \\
Electrical and Computer Engineering & Advanced Network Science Initiative \\ 
University of Washington, WA, USA  & Los Alamos National Laboratory, NM, USA \\
& \\
Christopher DeMarco & Ruisheng Diao  \\
Electrical and Computer Engineering & AI \& System Analytics \\ 
University of Wisconsin-Madison, WI, USA & GEIRI North America, CA, USA\\
& \\
Michael Ferris & St\'{e}phane Fliscounakis \\
Computer Sciences & Research and Development Division \\
University of Wisconsin-Madison, WI, USA & Reseau de Transport d'Electricite (RTE), Paris, France \\
& \\
Scott Greene & Renke Huang \\
Electrical and Computer Engineering & Electricity Infrastructure Group \\
University of Wisconsin-Madison, WI, USA & Pacific Northwest National Laboratory, WA, USA \\
& \\
 C\'{e}dric Josz  & Roman Korab \\
 Industrial Engineering and Operations Research & Electrical Engineering \\
 Columbia University, NY, USA & Silesian University of Technology, Gliwice, Poland \\
& \\
Bernard Lesieutre & Jean Maeght\\
Electrical and Computer Engineering & Research and Development Division \\ 
University of Wisconsin-Madison, WI, USA & Reseau de Transport d'Electricite (RTE), Paris, France \\
& \\
Terrence W. K. Mak & Daniel K. Molzahn \\
Industrial and Systems Engineering & Electrical and Computer Engineering \\
Georgia Institute of Technology, GA, USA & Georgia Institute of Technology, GA, USA \\
& \\
Thomas J. Overbye & Patrick Panciatici \\
Electrical and Computer Engineering & Research and Development Division \\
Texas A\&M University, TX, USA & Reseau de Transport d'Electricite (RTE), Paris, France \\
& \\
Byungkwon Park & Jonathan Snodgrass \\
Computational Sciences and Engineering Division & Electrical and Computer Engineering \\
Oak Ridge National Laboratory, TN, USA & University of Wisconsin-Madison, WI, USA \\
& \\
Ahmad Tbaileh & Pascal Van Hentenryck \\
Electricity Infrastructure Group & Industrial and Systems Engineering \\
Pacific Northwest National Laboratory, WA, USA & Georgia Institute of Technology, GA, USA\\
& \\
Ray Zimmerman & \\
Applied Economics and Management & \\
Cornell University, Ithaca, NY, USA & \\
\end{tabular}
\vspace{0.25cm}

The PGLib-OPF benchmarks are available at \url{https://github.com/power-grid-lib/pglib-opf.}\\
{\footnotesize $^{*}$ Corresponding author carleton@coffrin.com.}
}
}

\maketitle
\clearpage

\begin{abstract}
In recent years, the power systems research community has seen an explosion of novel methods for formulating the AC power flow equations. Consequently, benchmarking studies using the seminal AC Optimal Power Flow (AC-OPF) problem have emerged as the primary method for evaluating these emerging methods.  However, it is often difficult to directly compare these studies due to subtle differences in the AC-OPF problem formulation as well as the network, generation, and loading data that are used for evaluation.  To help address these challenges, this IEEE PES Task Force report proposes a standardized AC-OPF mathematical formulation and the PGLib-OPF networks for benchmarking AC-OPF algorithms.  A motivating study demonstrates some limitations of the established network datasets in the context of benchmarking AC-OPF algorithms and a validation study demonstrates the efficacy of using the PGLib-OPF networks for this purpose.  In the interest of scientific discourse and future additions, the PGLib-OPF benchmark library is open-access and all the of network data is provided under a creative commons license.
\end{abstract}

\begin{IEEEkeywords}
Nonlinear Optimization, Convex Optimization, AC Optimal Power Flow, Benchmarking
\end{IEEEkeywords}

\section*{Nomenclature}

\begin{IEEEdescription}[\IEEEusemathlabelsep\IEEEsetlabelwidth{$Y^s = g^s + \bm i$}]
  \item [{$N$}]  - The set of buses in the network 
  \item [{$G$}]  - The set of generators in the network 
  \item [{$E$}]  - The set of {\em from} branches in the network 
  \item [{$E^R$}]  - The set of {\em to} branches in the network 
  \item [{$\bm i$}] - Imaginary number constant
  \item [{$\bm e$}] - Exponential constant
  \item [{$S = p+ \bm iq$}] - AC power
  \item [{$V = v \angle \theta$}]  - AC voltage
  \item [{$Z = r + \bm i x$}] - Branch impedance
  \item [{$Y = g + \bm i b$}]  - Branch admittance
  \item [{$T = t \angle \theta^t$}]  - Branch transformer properties
  \item [{$Y^s = g^s + \bm i b^s$}]  - Bus shunt admittance
  \item [{$\dot{v}$}] - Nominal base voltage
  \item [{$b^c$}] - Line charging
  \item [{$s^u$}] - Branch apparent power limit
  \item [{$I^u$}] - Branch current magnitude limit
  \item [{$\theta^\Delta$}] - Voltage angle difference limit
  \item [{$S^d = p^d+ \bm iq^d$}] - AC power demand
  \item [{$S^g = p^g+ \bm iq^g$}] - AC power generation
  \item [{$c_0,c_1,c_2$}] - Generation cost coefficients 
   \item [{$\Re(\cdot), \Im(\cdot)$}] - Real and imag. parts of a complex number
   \item [{$|\cdot | \angle \cdot$}] - Magnitude and angle of a complex number
   \item [{$(\cdot)^*$}] - Conjugate of a complex number
  \item [{$x^l, x^u$}] - Lower and upper bounds of $x$, respectively
  \item [{$\bm x$}] - A constant value
  \item [{$\hat{x}$}] - An estimation of $x$
  \item [{$\mu$}] - Mean of a normal distribution
  \item [{$\sigma$}] - Standard deviation of a normal distribution
  \item [{$\lambda$}] - Rate of an exponential distribution
\end{IEEEdescription}

\section{Introduction}

Over the last decade, power systems research has experienced an explosion in variations of the steady-state AC power flow equations.  These include approximations such as the LPAC \cite{LPAC_ijoc}, IV-Flow \cite{IVModel} and relaxations such as the Second-Order Cone (SOC) \cite{1664986}, Convex-DistFlow (CDF) \cite{6102366}, Quadratic Convex (QC) \cite{qc_opf_tps,Hijazi2017}, Semidefinite Programming (SDP) \cite{Bai2008383}, and Moment / Sum-of-Squares Hierarchies \cite{7038397,7024950,6980142,moment_hierarchy}, just to name a few. Surveys of the power flow relaxation and approximation literature are provided in~\cite{coffrin_roald_videos,wang2018,molzahn_hiskens-fnt2017}. Much of the excitement underlying this line of research was ignited when \cite{5971792} demonstrated that the SDP relaxation could provide globally optimal solutions to a variety of the transmission system networks distributed with \matpower~\cite{matpower}.  Combining these results with industrial-strength convex optimization tools (e.g., Gurobi \cite{gurobi}, Cplex \cite{cplex}, Mosek \cite{mosek}) promises efficient and reliable algorithms for a wide variety of applications in power systems such as 
% economic dispatch \cite{ac_opf_origin}, 
optimal power flow \cite{ac_opf_origin,744492,744495,ferc4}, optimal transmission switching \cite{4492805}, and network expansion planning \cite{6407493}, just to name a few. 

Independent of the specific problem domain or the power flow model under consideration, all of these novel methods require AC power network data for experimental validation and it is important that suitable data are used to validate these emerging techniques \cite{hans_comp}.  However, due to the sensitive nature of critical infrastructure, detailed real-world network data is often difficult to obtain, even under non-disclosure agreements.  Consequently, many of the available power network datasets are over thirty years old (e.g. \cite{IEEEBench}), were originally designed for testing AC Power Flow algorithms, and lack the parameters needed for testing the AC Optimal Power Flow algorithms, e.g. branch thermal limits and generator cost functions.  The lack of comprehensive and modern AC power network data has been recognized by the Advanced Research Projects Agency-Energy (ARPA-e) and the GRID DATA program \cite{GRIDDATA} has resulted in a number of new network datasets, which are synthetically generated to match the statistics of real-world networks and provided as open-access.

To help improve the evaluation of AC-OPF algorithms, this IEEE PES Task Force report introduces the Power Grid Library for benchmarking the AC-OPF problem, \pglibopf.  \pglibopf~is a comprehensive collection of creative-commons AC transmission system networks curated in the \matpower~data format with all of the data required for modeling the proposed AC-OPF problem. This report is a companion document to the \pglibopf{} data.  Section \ref{sec:acopf} begins with a detailed specification of the AC-OPF problem.  Section~\ref{sec:motivation} provides a brief overview of the established \mbox{AC-OPF} cases provided with \matpower~and highlights some limitations in using these cases for benchmarking AC-OPF algorithms.  Section~\ref{sec:networks} provides a survey of known creative-commons transmission system network data and highlights the missing information in each network data source.  Section \ref{sec:data_models} briefly introduces methods for addressing the missing information in these datasets.  Section \ref{sec:pglibopf} introduces the \pglibopf~networks and conducts a baseline validation study to demonstrate that they are suitable for benchmarking the  AC-OPF problem. Section~\ref{sec:conclusion} provides concluding remarks. An appendix summarizes some variants of the proposed AC-OPF model and discusses the additional data that may be necessary in order to make the \pglibopf{} test cases applicable to other classes of power system optimization and control problems.

\section{The AC Optimal Power Flow Formulation}
\label{sec:acopf}

Many variations of the AC-OPF problem are relevant to power system analysis.  However, in accordance with \pglib~repository requirements, this section nominates a specific version of the AC-OPF problem for algorithmic benchmarking.  The following mathematical model presents a variant of the \mbox{AC-OPF} problem that is often used in related AC-OPF publications (e.g. \cite{1664986,6102366,qc_opf_tps,5971792,moment_hierarchy,LPAC_ijoc,IVModel}) and is readily encoded in the \mbox{\matpower}~network data format.

\begin{figure}[t]
\centering
\begin{circuitikz}[scale=1]\draw
  (0,2) node[anchor=south] {}
  (0,2) to[short,o-] (1.2,2)
  		to[cute inductor] (1.2,0)
  		to[short,-o](0,0)
  (0,2) to[open] (0,0)
  (3,0) to[short] (2,0)
  		to[cute inductor] (2,2)
  		to[short] (3,2)
  (1.65,2) node[anchor=south]{${\bm t_l \angle \bm \theta^t_l}\!:\!1$}
  (2,0) node[anchor=east] {}
        to[short] (7,0)
  (2,2) node[anchor=east] {}
  		to[short] (3,2)
  		to[R=$\bm r_l$] (5,2)
  		to[cute inductor,l=$\bm x_l$] (7,2)
  		to[C,l_=$\frac{\bm b^c_l}{2}$] (7,0)
  (3,0) to[C,l_=$\frac{\bm b^c_l}{2}$] (3,2)
  (7,2) to[short,-o] (8,2)
  (8,2) node[anchor=south] {}
  (7,0) to[short,-o] (8,0)
  (8,2) to[open] (8,0);
\end{circuitikz}
\caption{$\Pi$-circuit branch model with an ideal transformer. This is the branch model used by \matpower~\cite{matpower}.}
\label{fig:linemodel}
\end{figure}
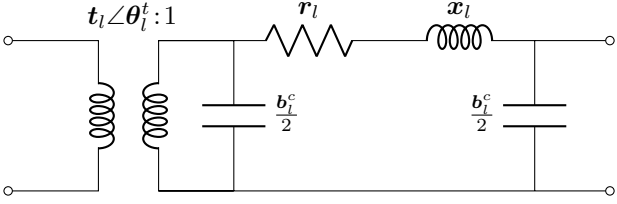

The proposed AC-OPF problem requires the following network parameters: a set of bus ids $N$; a set of branch ids $E$ with an arbitrary orientation; a set $E^R$ that captures the reverse orientation of the branches; and a set of generator ids $G$.  For each bus $i \in N$: the sets $E_i$ and $E^R_i$ indicate the subset of edges that are incident to that bus; the set $G_i$ reflects the subset of generator ids that are connected to that bus; $\bm {S^d}_i$ is the constant power demand; $\bm {Y^s}_i$ is the bus shunt admittance; and $\bm {v^l}_i, \bm {v^u}_i$ indicate the operating range of the bus' voltage magnitude.  For each generator $k \in G$: $\bm {S^{gl}}_k, \bm {S^{gu}}_k$ indicate the generator's power injection range; and $\bm c_{2k}, \bm c_{1k}, \bm c_{0k}$ provide the coefficients of a quadratic active power cost function.  For each branch $(l,i,j) \in E$: $i$ and $j$ are the {\em from} and {\em to} buses respectively and $l$ is the branch id; the series admittance, line charge, and transformer parameters are given by $\bm Y_{l}, \bm {b^c}_{l}, \bm T_{l}$ respectively; the branch's thermal limit in given by $\bm {s^u}_{l}$; and the branch voltage angle difference range is $\bm {\theta^{\Delta l}}_{l}, \bm {\theta^{\Delta u}}_{l}$.   Lastly, a voltage angle reference bus $\bm {ref} \in N$ is specified.  Most of these parameters are specified directly in a \matpower~data file.  However, the following parameters need to be computed from the raw data as follows,
\begin{subequations}
\begin{align}
    & \bm Y_{l} = \bm Z^{-1}_{l} = \frac{\bm r_l}{\bm r^2_l + \bm x^2_l} - \bm i \frac{\bm x_l}{\bm r^2_l + \bm x^2_l} \\
    & \bm T_l = \bm t_l \cos(\bm \theta^t_l) + \bm i \, \bm t_l \sin(\bm \theta^t_l)
\end{align}
\end{subequations}
% where $\bm r_l$, $\bm x_l$, $\bm \theta^t$, and $\bm t_l$  denote the series resistance, series reactance, phase shift, and tap ratio, respectively, of branch $l$.
Figure~\ref{fig:linemodel} shows the circuit model used to represent each branch.

Model \ref{model:acopf} presents the AC-OPF problem as a non-convex nonlinear mathematical program over complex values and variables. A detailed description of the model's notation and derivation can be found in \cite{qc_opf_tps}.  The objective function \eqref{eq:acopf:obj} strives to minimize the cost of active power injections.  Constraint \eqref{eq:acopf:1} fixes the voltage angle of the reference bus.  Constraint \eqref{eq:acopf:2} sets the generator injection limits and constraint \eqref{eq:acopf:3} sets the bus voltage magnitude limits.  Constraint \eqref{eq:acopf:4} captures the nodal power balance and constraints \eqref{eq:acopf:5}--\eqref{eq:acopf:6} ensure that the branch power flows are consistent with Ohm's Law.  Finally, constraints \eqref{eq:acopf:7} and \eqref{eq:acopf:8} capture the branch thermal and voltage angle difference limits.  It is important to note that solving Model \ref{model:acopf} is NP-Hard \cite{verma2009power} in general, even if the network has a tree topology \cite{7063278}.  Consequently, it is expected that solution methods for  Model \ref{model:acopf} will exhibit a wide variety of quality-runtime tradeoffs and will be specialized to different classes of inputs.

\begin{model}[t]
\caption{The AC Optimal Power Flow Problem (AC-OPF)}
\label{model:acopf}
\vspace{-0.3cm}
\begin{subequations}
\begin{align}
&\mbox{\bf variables: } \nonumber \\
& S^g_k \;\; \forall k\in G \nonumber \\
& V_i \;\; \forall i\in N \nonumber \\
& S_{lij} \;\; \forall (l,i,j) \in E \cup E^R \nonumber \\
&\mbox{\bf minimize: } \sum_{k \in G} \bm c_{2k} (\Re(S^g_k))^2 + \bm c_{1k}\Re(S^g_k) + \bm c_{0k} \label{eq:acopf:obj}\\
&\mbox{\bf subject to: } \nonumber \\
& \angle V_{\bm {ref}} = 0 \label{eq:acopf:1} \\
& \bm {S^{gl}}_k \leq S^g_k \leq \bm {S^{gu}}_k \;\; \forall k \in G \label{eq:acopf:2}  \\
& \bm {v^l}_i \leq |V_i| \leq \bm {v^u}_i \;\; \forall i \in N \label{eq:acopf:3} \\
& \sum_{\substack{k \in G_i}} S^g_k - \bm {S^d}_i - \bm {Y^s}_{i} |V_i|^2 = \sum_{\substack{(l,i,j)\in E_i \cup E_i^R}} S_{lij} \;\; \forall i\in N \label{eq:acopf:4} \\ 
& S_{lij} = \left( \bm Y^*_{l} - \bm i\frac{\bm {b^c}_{l}}{2} \right) \frac{|V_i|^2}{|\bm{T}_{l}|^2} - \bm Y^*_{l} \frac{V_i V^*_j}{\bm{T}_{l}} \;\; \forall (l,i,j)\in E \label{eq:acopf:5} \\
& S_{lji} = \left( \bm Y^*_{l} - \bm i\frac{\bm {b^c}_{l}}{2} \right) |V_j|^2 - \bm Y^*_{l} \frac{V^*_i V_j}{\bm{T}^*_{l}} \;\; \forall (l,i,j)\in E \label{eq:acopf:6} \\
& |S_{lij}| \leq \bm {s^u}_{l} \;\; \forall (l,i,j) \in E \cup E^R  \label{eq:acopf:7} \\
& \bm {\theta^{\Delta l}}_{l} \leq \angle (V_i V^*_j) \leq \bm {\theta^{\Delta u}}_{l} \;\; \forall (l,i,j) \in E  \label{eq:acopf:8}
\end{align}
\end{subequations}
\end{model}

%
% Table generated with: Julia v1.0.2, JuMP v0.18.5, PowerModels v0.9.3, Ipopt v0.5.1 and HSL MA27
% Repo: matpower-experiments
% Result: pglib-tbl.py
%

\begin{table*}[h]
\center
\caption{AC-OPF Optimality Gaps on Network Datasets Distributed with \matpower~v6.0.}
\begin{tabular}{|r|r|r||r||r||r|r|r|r|r|r|r|r|r|r|r|r|}
\hline
& & & \$/h & Gap (\%) & \multicolumn{2}{c|}{Runtime (seconds)} \\
Test Case & $|N|$ & $|E|$ & AC & SOC & AC & SOC \\
\hline
\hline
case5 & 5 & 6 & 1.7552e+04 & 14.55 & $<$1 & $<$1 \\
\hline
case6ww & 6 & 11 & 3.1440e+03 & 0.63 & $<$1 & $<$1 \\
\hline
case9 & 9 & 9 & 5.2967e+03 & 0.01 & $<$1 & $<$1 \\
\hline
case9target & 9 & 9 & n.s. & inf. & $<$1 & $<$1 \\
\hline
case14 & 14 & 20 & 8.0815e+03 & 0.08 & $<$1 & $<$1 \\
\hline
case24\_ieee\_rts & 24 & 38 & 6.3352e+04 & 0.02 & $<$1 & $<$1 \\
\hline
case30 & 30 & 41 & 5.7689e+02 & 0.58 & $<$1 & $<$1 \\
\hline
case\_ieee30 & 30 & 41 & 8.9061e+03 & 0.05 & $<$1 & $<$1 \\
\hline
case39 & 39 & 46 & 4.1864e+04 & 0.03 & $<$1 & $<$1 \\
\hline
case57 & 57 & 80 & 4.1738e+04 & 0.07 & $<$1 & $<$1 \\
\hline
case89pegase & 89 & 210 & 5.8198e+03 & 0.17 & $<$1 & $<$1 \\
\hline
case118 & 118 & 186 & 1.2966e+05 & 0.25 & $<$1 & $<$1 \\
\hline
case145 & 145 & 453 & n.s. &  inf. & 14 & 7 \\
\hline
case\_illinois200 & 200 & 245 & 3.6748e+04 & 0.02 & $<$1 & $<$1 \\
\hline
case300 & 300 & 411 & 7.1973e+05 & 0.15 & $<$1 & $<$1 \\
\hline
case1354pegase & 1354 & 1991 & 7.4069e+04 & 0.08 & 4 & 5 \\
\hline
%case1888rte & 1888 & 2531 & 5.9805e+04 & 0.38 & 3349 & 61 \\
%\hline
case1951rte & 1951 & 2596 & 8.1738e+04 & 0.08 & 17 & 26 \\
\hline
case2383wp & 2383 & 2896 & 1.8685e+06 & 1.05 & 9 & 6 \\
\hline
case2736sp & 2736 & 3504 & 1.3079e+06 & 0.30 & 8 & 5 \\
\hline
case2737sop & 2737 & 3506 & 7.7763e+05 & 0.26 & 6 & 4 \\
\hline
case2746wop & 2746 & 3514 & 1.2083e+06 & 0.37 & 7 & 4 \\
\hline
case2746wp & 2746 & 3514 & 1.6318e+06 & 0.33 & 7 & 5 \\
\hline
case2848rte & 2848 & 3776 & 5.3022e+04 & 0.08 & 46 & 7 \\
\hline
case2868rte & 2868 & 3808 & 7.9795e+04 & 0.07 & 28 & 8 \\
\hline
case2869pegase & 2869 & 4582 & 1.3400e+05 & 0.09 & 9 & 58 \\
\hline
case3012wp & 3012 & 3572 & 2.5917e+06 & 0.78 & 12 & 6 \\
\hline
case3120sp & 3120 & 3693 & 2.1427e+06 & 0.54 & 11 & 6 \\
\hline
case3375wp & 3374 & 4161 & 7.4120e+06 & 0.26 & 13 & 98 \\
\hline
case6468rte & 6468 & 9000 & 8.6829e+04 & 0.23 & 67 & 226 \\
\hline
case6470rte & 6470 & 9005 & 9.8345e+04 & 0.17 & 61 & 105 \\
\hline
case6495rte & 6495 & 9019 & 1.0628e+05 & 0.45 & 44 & 239 \\
\hline
case6515rte & 6515 & 9037 & 1.0980e+05 & 0.38 & 45 & 36 \\
\hline
case9241pegase & 9241 & 16049 & 3.1591e+05 & 1.75 & 61 & 230 \\
\hline
case13659pegase & 13659 & 20467 & 3.8611e+05 & 1.52 & 228 & 215 \\
\hline
\end{tabular}
\label{tbl:gaps_time_mp}
\end{table*}
\noindent

\begin{table*}[!ph]
\center
\footnotesize
%\small
\vspace{-0.5cm}
\caption{A Survey of Transmission System Data Sources$^*$}
\begin{tabular}{|r||c|c|c|c|c|c|c|r|r|r|r|r|r|r|c|c|}
\hline
 & Original  & Generator  & Generator  & Thermal  \\
Name &  Source &  Injection Limits &  Costs &  Limits \\
\hline
%%%%%%%%%%%%%%%%%%
\hline
\multicolumn{5}{|c|}{Publication Test Cases}\\
\hline
3-Bus& \cite{6120344} & \cite{6120344}  & \cite{6120344}  & \cite{6120344}  \\
\hline
case5 & \cite{5589973} & \cite{5589973} & \cite{5589973}  & --- \\
\hline
case30-as & \cite{4075418} & \cite{4075418} & \cite{4075418} & \cite{4075418} \\
\hline
%case30-fsr & \cite{4075418} & \cite{630479} & \cite{630479} & \cite{4075418} \\
%\hline
case39 & \cite{Bills_1970} & \cite{Bills_1970, pai1989energy} & \cite{4562306} & \cite{TCCalculator} \\
\hline
\hline
\multicolumn{5}{|c|}{IEEE Power Flow Test Cases}\\
\hline
14 Bus & \cite{IEEEBench} & --- & --- & --- \\
\hline
30 Bus & \cite{IEEEBench} & --- & --- & --- \\
\hline
57 Bus & \cite{IEEEBench} & --- & --- & --- \\
\hline
118 Bus & \cite{IEEEBench} & --- & --- & --- \\
\hline
300 Bus & \cite{IEEEBench} & --- & --- & --- \\
\hline
\hline
\multicolumn{5}{|c|}{IEEE Dynamic Test Cases}\\
\hline
17 Generator & \cite{IEEEBench} & --- & --- & --- \\
\hline
\hline
\multicolumn{5}{|c|}{IEEE Reliability Test Systems (RTS)}\\
\hline
RTS-79 & \cite{4113721} & \cite{4113721} & \cite{4113721,rts_79_gen_data} & \cite{4113721} \\
\hline
RTS-96 & \cite{rts96} & \cite{4113721} & \cite{rts_79_gen_data} & \cite{4113721} \\
\hline
\hline
\multicolumn{5}{|c|}{Polish Test Cases}\\
\hline
case2383wp & \cite{matpower} & \cite{matpower} & \cite{matpower} & \cite{matpower} \\
\hline
case2736sp & \cite{matpower} & \cite{matpower} & \cite{matpower} & \cite{matpower} \\
\hline
case2737sop & \cite{matpower} & \cite{matpower} & \cite{matpower} & \cite{matpower} \\
\hline
case2746wop & \cite{matpower} & \cite{matpower} & \cite{matpower} & \cite{matpower} \\
\hline
case2746wp & \cite{matpower} & \cite{matpower} & \cite{matpower} & \cite{matpower} \\
\hline
case3012wp & \cite{matpower} & \cite{matpower} & \cite{matpower} & --- \\
\hline
case3120sp & \cite{matpower} & \cite{matpower} & \cite{matpower} & ---\\
\hline
case3375wp & \cite{matpower} & --- & \cite{matpower} & --- \\
\hline
\hline
\multicolumn{5}{|c|}{PEGASE Test Cases}\\
\hline
case89pegase & \cite{rte_cases} & \cite{rte_cases} & --- &  \cite{rte_cases}, partial \\
\hline
case1354pegase & \cite{rte_cases} & \cite{rte_cases} & --- & \cite{rte_cases}, partial \\
\hline
case2869pegase & \cite{rte_cases} & \cite{rte_cases} & --- & \cite{rte_cases}, partial \\
\hline
case9241pegase & \cite{rte_cases} & \cite{rte_cases} & --- & \cite{rte_cases}, partial \\
\hline
case13659pegase & \cite{rte_cases} & \cite{rte_cases} & --- & \cite{rte_cases}, partial \\
\hline
\hline
\multicolumn{5}{|c|}{RTE Test Cases}\\
\hline
case1888rte & \cite{rte_cases} & \cite{rte_cases} & --- &  \cite{rte_cases}, partial \\
\hline
case1951rte & \cite{rte_cases} & \cite{rte_cases} & --- & \cite{rte_cases}, partial \\
\hline
case2848rte & \cite{rte_cases} & \cite{rte_cases} & --- & \cite{rte_cases}, partial \\
\hline
case2868rte & \cite{rte_cases} & \cite{rte_cases} & --- & \cite{rte_cases}, partial \\
\hline
case6468rte & \cite{rte_cases} & \cite{rte_cases} & --- & \cite{rte_cases}, partial \\
\hline
case6470rte & \cite{rte_cases} & \cite{rte_cases} & --- & \cite{rte_cases}, partial \\
\hline
case6495rte & \cite{rte_cases} & \cite{rte_cases} & --- & \cite{rte_cases}, partial \\
\hline
case6515rte & \cite{rte_cases} & \cite{rte_cases} & --- & \cite{rte_cases}, partial \\
\hline
\hline
\multicolumn{5}{|c|}{ACTIVSg Test Cases}\\
\hline
ACTIVSg200 & \cite{7725528} & \cite{7725528} & \cite{7725528} & \cite{7725528} \\
\hline
% ACTIVSg500 &  \cite{7725528} & \cite{7725528} & \cite{7725528} & \cite{7725528}\\
% \hline
% ACTIVSg2000 &  \cite{7725528} & \cite{7725528} & \cite{7725528} & \cite{7725528} \\
% \hline
% ACTIVSg10k &  \cite{7725528} & \cite{7725528} & \cite{7725528} & \cite{7725528} \\
% \hline
%
\hline
\multicolumn{5}{|c|}{Sustainable Data Evolution Technology Test Cases}\\
\hline
SDET  500 & \cite{SDETBench} & \cite{SDETBench} & --- & \cite{SDETBench} \\
\hline
%SDET 2000 & \cite{SDETBench} & \cite{SDETBench} & --- & \cite{SDETBench} \\
%\hline
SDET 3000 & \cite{SDETBench} & \cite{SDETBench} & --- & \cite{SDETBench} \\
\hline
SDET 4000 & \cite{SDETBench} & \cite{SDETBench} & --- & \cite{SDETBench} \\
\hline
\hline
\multicolumn{5}{|c|}{Power Systems Engineering Research Center Test Cases}\\
\hline
WECC 240 Bus & \cite{6039476,PSERCProjects} & \cite{PSERCProjects} & \cite{PSERCProjects} & \cite{PSERCProjects} \\
\hline
\hline
\multicolumn{5}{|c|}{Grid Optimization Competition Test Cases}\\
\hline
179 Bus & \cite{GOCBench} & \cite{GOCBench} & --- & --- \\
\hline
Network\_02, 173 scenarios & \cite{GOCBench} & \cite{GOCBench} & \cite{GOCBench} & \cite{GOCBench} \\
\hline
Network\_03, 200 scenarios & \cite{GOCBench} & \cite{GOCBench} & \cite{GOCBench} & \cite{GOCBench} \\
\hline
Network\_06, 124 scenarios & \cite{GOCBench} & \cite{GOCBench} & \cite{GOCBench} & \cite{GOCBench} \\
\hline
Network\_08, 373 scenarios & \cite{GOCBench} & \cite{GOCBench} & \cite{GOCBench} & \cite{GOCBench} \\
\hline
Network\_09, 195 scenarios & \cite{GOCBench} & \cite{GOCBench} & \cite{GOCBench} & \cite{GOCBench} \\
\hline
Network\_12, 50 scenarios & \cite{GOCBench} & \cite{GOCBench} & \cite{GOCBench} & \cite{GOCBench} \\
\hline
Network\_13, 103 scenarios & \cite{GOCBench} & \cite{GOCBench} & \cite{GOCBench} & \cite{GOCBench} \\
\hline
Network\_14, 50 scenarios & \cite{GOCBench} & \cite{GOCBench} & \cite{GOCBench} & \cite{GOCBench} \\
\hline
Network\_20, 100 scenarios & \cite{GOCBench} & \cite{GOCBench} & \cite{GOCBench} & \cite{GOCBench} \\
\hline
Network\_25, 60 scenarios & \cite{GOCBench} & \cite{GOCBench} & \cite{GOCBench} & \cite{GOCBench} \\
\hline
Network\_30, 48 scenarios & \cite{GOCBench} & \cite{GOCBench} & \cite{GOCBench} & \cite{GOCBench} \\
\hline
Network\_70, 422 scenarios & \cite{GOCBench} & \cite{GOCBench} & \cite{GOCBench} & \cite{GOCBench} \\
\hline
Network\_75, 40 scenarios & \cite{GOCBench} & \cite{GOCBench} & \cite{GOCBench} & \cite{GOCBench} \\
\hline
Network\_82, 10 scenarios & \cite{GOCBench} & \cite{GOCBench} & \cite{GOCBench} & \cite{GOCBench} \\
\hline
Network\_82, 40 scenarios & \cite{GOCBench} & \cite{GOCBench} & \cite{GOCBench} & \cite{GOCBench} \\
\hline
Network\_83, 130 scenarios & \cite{GOCBench} & \cite{GOCBench} & \cite{GOCBench} & \cite{GOCBench} \\
\hline
Network\_86, 40 scenarios & \cite{GOCBench} & \cite{GOCBench} & \cite{GOCBench} & \cite{GOCBench} \\
\hline
Network\_88, 30 scenarios & \cite{GOCBench} & \cite{GOCBench} & \cite{GOCBench} & \cite{GOCBench} \\
\hline
\end{tabular}
\\
\vspace{0.1cm}
$*$ - only creative commons data sources were considered
\label{tbl:opf:inst}
\end{table*}
%\clearpage

\section{Motivation}
\label{sec:motivation}

To motivate the need for a careful curation of the network data in \pglibopf, this section conducts a preliminary study of thirty-five AC transmission system datasets that are distributed with \matpower~v6.0 \cite{matpower,zimmerman_ray_d_2016_3237810}.  The {\em optimality gap} measure is used as a simple and preliminary test of AC-OPF difficulty, as one expects that challenging cases will exhibit a large optimality gap.  Given a feasible solution to the AC-OPF problem and the solution to a convex relaxation, the optimality gap is defined as the relative difference between the objective values of the feasible solution and the relaxation: %, as defined in~\eqref{eq:opt_gap}.
\begin{align}
\frac{\text{AC Heuristic} - \text{AC Relaxation}}{\text{AC Heuristic}} \label{eq:opt_gap}
\end{align}
There are a wide variety of both AC heuristics (i.e., methods for obtaining feasible solutions to AC OPF problems~\cite{744492,744495,ferc4}) and convex relaxation techniques~\cite{coffrin_roald_videos,wang2018,molzahn_hiskens-fnt2017}. In the interest of simplicity, this preliminary study will use a nonlinear optimization solver that converges to a KKT point as a heuristic for finding AC feasible solutions and a simple Second-Order Cone (SOC) relaxation \cite{1664986} for providing objective bounds.  All of the results were computed using IPOPT 3.12 \cite{Ipopt} with the HSL \cite{hsl_lib} linear algebra library on a server with two 2.10GHz Intel CPU and  128GB of RAM.  PowerModels.jl v0.9 \cite{8442948} was used to formulate and solve both mathematical programs.

The results of this study are presented in Table \ref{tbl:gaps_time_mp}.  The data highlights two core points: (1) By-in-large the optimality gaps are less than 1\%. Although a large optimality gap is not a necessary condition for AC-OPF hardness, it provides a good indication of a challenging instance.  This work will demonstrate that much more significant gaps are possible, providing a significant increase in the variety of network cases for AC-OPF algorithm benchmarking; (2) No feasible solution was found in two cases, case9target and case145.  This could suggest that these cases are challenging for AC heuristics.  However, the SOC relaxation provides a numerical proof that these cases have no feasible AC-OPF solution,\footnote{That is, there does not exist an assignment of the variables that can simultaneously satisfy all of the constraints in Model \ref{model:acopf}.} suggesting that data quality is the source of infeasibility and not algorithmic difficulty.
%(3) The AC heuristic runtime of case1888rte is a notable outlier, which may be numerically challenging.  However, the subsequent studies in this report will demonstrate this is an anomaly that is eliminated by improved data quality considerations.  
%
Overall, this simplistic study demonstrates some of the shortcomings of focusing exclusively on the network data that is distributed with \matpower~v6.0 for benchmarking AC-OPF algorithms.  The careful curation of AC network data developed in the following sections results in modified network datasets featuring significant optimality gaps, which will help to emphasize the differences between various AC-OPF solution methods.

\section{Publicly Available Network Data}
\label{sec:networks}

In the interest of curating a comprehensive collection of \mbox{AC-OPF} networks, we begin with a survey of available network datasets.  To the best of our knowledge, Table \ref{tbl:opf:inst} summarizes all of the readily available transmission system datasets.\footnote{Only creative commons datasets were considered to comply with \pglib~data requirements.}  A careful investigation of these datasets reveals that very few networks include all of the data required to study Model~\ref{model:acopf}.  Table \ref{tbl:opf:inst} highlights the source of, or lack of, generation capacity limits, generation cost functions, and branch thermal limits in these datasets.  Cells containing ``---" indicate missing data that must be added before the network will be suitable for benchmarking Model \ref{model:acopf}.   To address the information that is missing in these datasets, the next section reviews a number of data-driven models that can be used to fill in these gaps.

\subsection{Network Omissions}

Some notable networks have not been included in Table \ref{tbl:opf:inst} for the following reasons.

\subsubsection{IEEE 30 Bus ``New England'' Dynamic Test System}
This test case is nearly identical to the IEEE 30 test case and would not bring additional value to the proposed collection of cases.

\subsubsection{30 Bus System of Ferrero, Shahidehpour, Ramesh}
In \pglibopf~ v20.07 the  30 Bus System of Ferrero, Shahidehpour and Ramesh \cite{630479} was removed due to being very similar to the 30 Bus System of Alsac and Stott \cite{4075418} and would not bring additional value to the proposed collection of cases.

\subsubsection{IEEE 50 Generator Dynamic Test System}
In its specified state, this test case does not converge to an AC power flow solution.  However, if the active generation upper bounds are increased to 1.5 times their given value and voltage bounds are set to $1 \pm 0.16$, then a solution can be obtained.  This solution still exhibits significant voltage drops, atypical of other networks.  Many of the lines in this network have negative $r$ and $x$ values, which is likely the result of a network reduction procedure. Also, the size of the generating units are one or two orders of magnitude larger than any documented generation unit in the U.S., which suggests that these ``generators'' may actually be modeling imports and exports of power. Since many of these characteristics are atypical compared to the other test cases, this network is omitted. 

\subsubsection{SDET 2000 System}
This system was omitted starting with \pglibopf~ v20.07, due to being similar to the Network 70 systems from the Grid Optimization Competition \cite{GOCBench}.

\subsubsection{ACTIVGs 500, 2000 and 10000 Bus Systems}
These three systems were omitted starting with \pglibopf~ v20.07, due to being similar to the Network 02, Network 06 and Network 13 systems from the Grid Optimization Competition \cite{GOCBench}.

%\subsection{Network Observations}
%
%\paragraph{IEEE 300 Test System}
%This test case exhibits very large voltage drops, with a voltage drop across branch 191--192 of 0.09 volts p.u.\ and 20 degrees. Such drops are not observed in the larger real-world networks \cite{Purchala:2005gt}, which brings the realism of this network into question.  Also note that in the \matpower~version of this network a phase shifting transformer is missing between buses 196 and 2040.
%
%\paragraph{IEEE RTS-96}
%Note that RTS-96 \cite{rts96} is simply three RTS-79 \cite{4113721} networks connected together.  One of the subtle points in these networks is the line charging on line $106$--$110$, which is 2.459 p.u\@.  In RTS-79, a -1.0 p.u.\ bus-shunt is used at bus 6 to balance out this massive line charge.  This value must be used in all three analogous parts of the RTS-96 network to produce a feasible AC power flow.  The generation costs presented in \cite{rts96,4113721} are not provided in \$/h, hence we adopt to interpretation of \cite{rts_79_gen_data} to get them in the correct units.
%
%\paragraph{3000 Bus Polish Test Cases}
%Unlike the other Polish network cases, these test cases have 10 lines with negative resistance values.  Hence, these lines can inject power into the network and invalidate some models, such as the AC copper plate relaxation.  However, these negative values do not appear to adversely effect the  power flow models considered in this report. 
%

\section{Data Driven Models}
\label{sec:data_models}

Ideally, the data missing from Table \ref{tbl:opf:inst} would be incorporated by returning to the original network design documents and extracting the required information, such as a generator's nameplate capacity and line conductor specifications.  Unfortunately, due to the age or the synthetic nature of these test cases, this approach to data completion is impractical.  This work proposes to leverage publicly available data sources to build data-driven models that can complete the missing data.  Such models may not reflect any specific real-world network, but at least they will reflect many of the statistical features found in realistic networks.  As identified in Section \ref{sec:networks}, the key pieces of missing data are generator injection limits, generation costs functions, and branch thermal limits.  The rest of this section reviews several data-driven models proposed in \cite{nesta}, which can be used to fill the gaps in these networks.

\subsection{Generator Models}

Most AC transmission datasets are brief in their description of the generation units.  Typically, only the active power injection limits, reactive power injection limits, and a specific generation dispatch point are provided.  A notable omission is an active power generation cost function, which is critical in formulating the objective in Model \ref{model:acopf}.  Two key observations can be used to address the limited information on generation units: (1) The U.S.\ Energy Information Administration (EIA) collects extensive data on generation units throughout the United States.  Two reports are particularly useful to this work, the detailed generator data (EIA-860 2012 \cite{EIA860}) and state fuel cost data (SEDS \cite{SEDS2012}); (2) The bulk of a generator's properties are driven by its mechanical design, which is in turn significantly influenced by its fuel type.  This work begins by developing a data-driven model for generators by assigning them a fuel category.  Once a fuel category is determined, probabilistic models for both fuel costs and power injection limits can be derived from publicly available data sources.  In the interest of simplicity, this work focuses on the four primary dispatchable fuel types, Petroleum (PEL), Natural Gas (NG), Coal (COW), and Nuclear Fuel (NUC).

\subsubsection{The Generation Fuel Category Model}

\begin{figure}[h!]
\center
    \includegraphics[width=9.5cm]{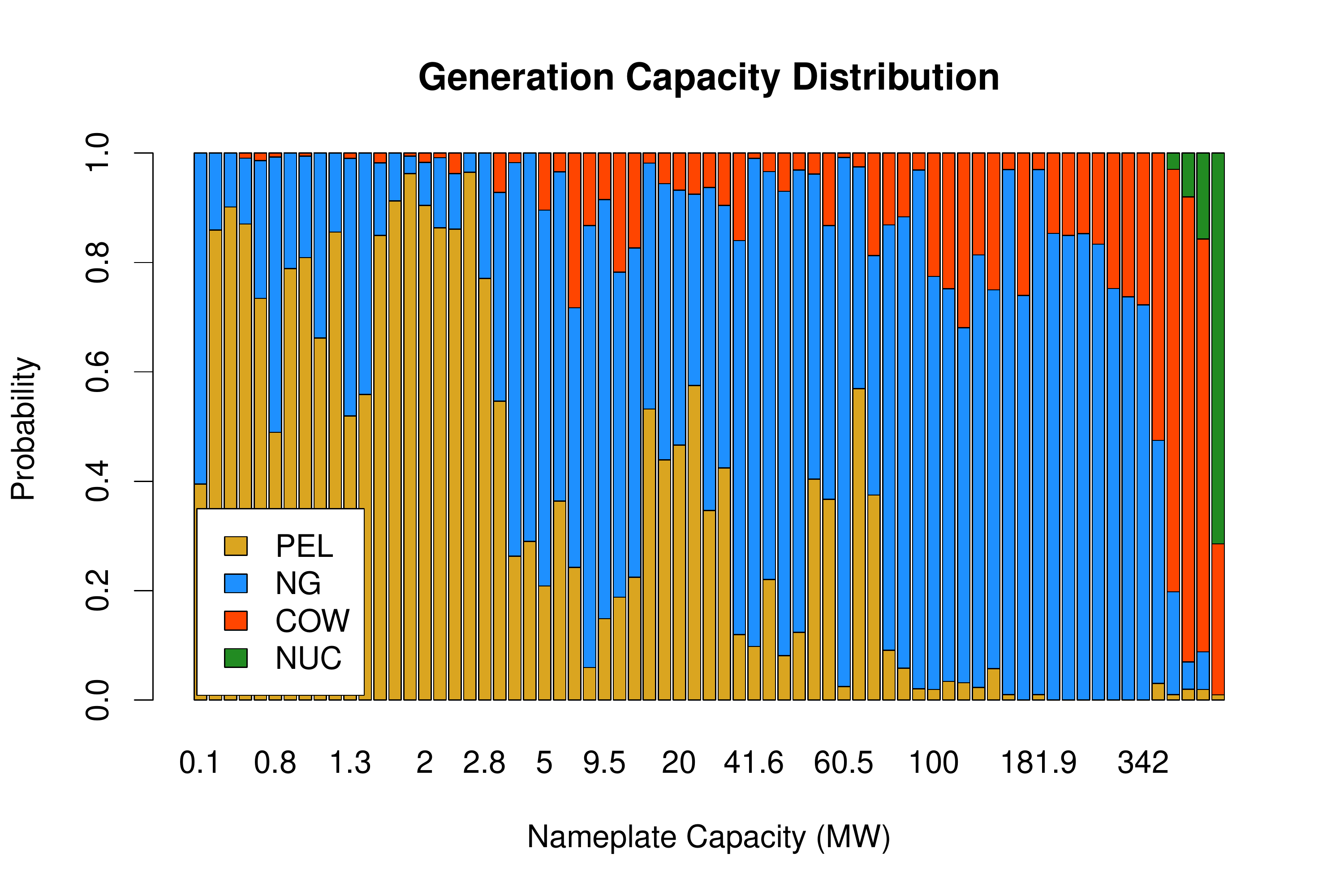}
    %\hspace{-0.4cm}
    \caption{An Empirical Distribution of Generation Fuel Categories by Nameplate Capacity.}
    %\vspace{-0.4cm}
\label{fig:gen_cap_bins}
\end{figure}

Assuming that a generator's active power injection limit (i.e. $\Re(\bm {S^{gu}})$) is a sufficient proxy for its nameplate capacity, one can use the empirical distribution presented in Figure \ref{fig:gen_cap_bins} for making a probabilistic guess of the fuel type of a given generation unit.  The fuel category classifier is built by selecting the corresponding nameplate capacity bin in Figure \ref{fig:gen_cap_bins} and rolling a weighted die to select the fuel type.
 
One important point is the identification of synchronous condensers.  In Model \ref{model:acopf} such devices are not explicitly identified but are modeled as generators with no active generation capabilities.  To identify such devices, one can introduce a new fuel category (e.g. SYNC), and any generator with active generation upper and lower bounds of 0 is assigned this category.  The empirical distribution shown in Figure \ref{fig:gen_cap_bins} combined with this special case forms the generation fuel classification model ({\bf GF-Stat}).

%There is also one special case to be considered in generator category classification.  Because many of these networks were originally designed as AC power flow test cases, the generators at the slack bus sometimes have no active generation value specified, as this would be a free variable in a power flow study.  If this occurs, we assume the generator is a very large and cheap power plant and assign it the class NUC\@.  

\subsubsection{Generation Capacity Models}

Some AC transmission system datasets lack reasonable generation injection limits, especially in cases that where originally designed for benchmarking AC Power Flow algorithms.  Two simple generation capability models are developed to address such cases.

\paragraph{Active Generation Capability}

An investigation of the nameplate capacities of each fuel type in \cite{nesta} revealed that an exponential distribution is a suitable model for PEL, NG, and COW generators, while a normal distribution is suitable for NUC generators.  The parameters from maximum likelihood estimation of these distributions are presented in Table \ref{tbl:gen:cap:model}.  Using these distributions the active generation capacity model ({\bf AG-Stat}) is constructed as follows.  Given a fuel category $f$ and an active generation upper limit or present output $\Re(\bm {S^g})$, the fuel category nameplate capacity distribution is sampled as $\bm {p^{gu}}$, until $\bm {p^{gu}} > \Re(\bm {S^{g}})$.  Then the maximum active power generation capacity is updated to the sampled value, i.e. $\Re(\bm {S^{gu}}) = \bm {p^{gu}}$.

\begin{table}[h!]
\center
%\footnotesize
\caption{Distribution Parameters for Generator Capacity Models.}
\begin{tabular}{|l||r|r|r|r|r|r|r|r|r|r|r|r|r|r|c|c|}
\hline
                   & Nameplate Capacity (MW) \\
 Fuel  Type & $\hat{\lambda}$ \\
\hline
\hline
%%%%%%%%%%%%%%%%%%%%%%%%%%%%%%%%%%%
PEL & 0.023254  \\
\hline
NG & 0.009188 \\
\hline
COW & 0.003201 \\
\hline
%NUC & 0.000957 \\
%\hline
%%%%%%%%%%%%%%%%%%%%%%%%%%%%%%%%%%%%
\end{tabular}

\vspace{0.2cm}
\begin{tabular}{|l||r|r|r|r|r||r|r|r|r|r|r|r|r|r|c|c|}
\hline
                   & \multicolumn{2}{c|}{Nameplate Capacity (MW)} \\
 Fuel  Type & $\hat{\mu}$ & $\hat{\sigma}$ \\
\hline
\hline
%%%%%%%%%%%%%%%%%%%%%%%%%%%%%%%%%%%
NUC & 1044.56 & 219.27 \\
\hline
%%%%%%%%%%%%%%%%%%%%%%%%%%%%%%%%%%%%
\end{tabular}
\label{tbl:gen:cap:model}
\end{table}

% \begin{figure}[h!]
% \center
%  \includegraphics[width=6.5cm]{results/gen_capacity_diff_rel_fuel_type_box_selected.pdf}
%       %\hspace{-0.4cm}
%     \caption{Distribution of Summer Peak Generation Capacity Reduction by Fuel Category.}
%     %\vspace{-0.4cm}
% \label{tbl:gen:cap:peak:dist}
% \end{figure}

\paragraph{Reactive Generation Capabilities}

In synchronous machines, reactive generation capabilities are tightly coupled with active generation capabilities.  Lacking detailed information about the generator's specifications it is observed in \cite{nesta} that the reactive power capability of a synchronous machine is roughly $\pm 50\%$ of its nameplate capacity leading to the \mbox{\bf RG-AM50} model.  This model assumes the given reactive power bounds are accurate, unless they exceed 50\% of the nameplate capacity, in which case, they are reduced to $\pm 50\%$ of the nameplate capacity.  This provides a pessimistic model of the generator's capabilities.

\subsubsection{Generation Cost Model}

Observing the simplicity of the cost function in Model \ref{model:acopf}, this work proposes to focus on the marginal costs of power generation in an idealized non-competitive environment.  Fuel price information is available in the ``Primary Energy, Electricity, and Total Energy Price Estimates, 2012" in the SEDS dataset \cite{SEDS2012}.  In \cite{nesta} it was observed that for the fuel categories of interest, the fuel costs are roughly normally distributed, with the parameters specified in Table \ref{tbl:gen:cost:model}.  Using these distributions the model for generation costs ({\bf AC-Stat}) is built as follows.  The fuel cost parameters from Table \ref{tbl:gen:cost:model} are assumed to be representative of the price variations across the various generating units, so that, given a fuel category, one simply draws a sample from the associated normal distribution to produce a linear fuel cost value (\$/BTU) for that generator. The conversion from heat energy input (BTU) to electrical energy output (MWh) depends on a generator's \emph{heat rate}, which differs based on the efficiency of the generating plant. Representative heat rate values are adopted from 2016 EIA data~\cite{eia_heatrate}. 

\begin{table}[h!]
\center
%\footnotesize
\caption{Generator Cost Model.}
\begin{tabular}{|r|r||r|r|r|r||r|r||r|r|r|r|r|r|r|c|c|}
\hline
& & \multicolumn{2}{c|}{Cost (\$/MWh)} \\
Fuel Category & SEDS Label & $\hat{\mu}$ & $\hat{\sigma}$ \\
\hline
\hline
%%%%%%%%%%%%%%%%%%%%%%%%%%%%%%%%%%%
PEL & Distillate Fuel Oil & 111.3398 & 9.6736 \\
\hline
NG & Natural Gas & 34.2731 & 10.9810 \\
\hline
COW & Coal & 24.7919 & 8.0866 \\
\hline
NUC & Nuclear Fuel & 7.2504 & 0.7534 \\
\hline
 %%%%%%%%%%%%%%%%%%%%%%%%%%%%%%%%%%%%
\end{tabular}
\label{tbl:gen:cost:model}
\end{table}

\subsection{Branch Thermal Limit Models}

Determining a transmission line's operational thermal rating is a intricate tasks that combines a wide variety of information such as the conductor's design, location, age and the season of the year.  Unfortunately, this information is unavailable in all of the network datasets presented in Table \ref{tbl:opf:inst}.  The only available branch parameters are the impedance ($Z = r+ \bm i x$ p.u.), line charge ($b^c$ p.u.), and often the nominal voltage ($\dot{v}$) on the connecting buses.  This work leverages two models for producing reasonable thermal limits, a data-driven approach and an arithmetic approach leveraging other model parameters.

\subsubsection{A Statistical Model}

After reviewing a number of network datasets with realistic thermal limits, \cite{nesta} concluded that the following exponential model was reasonable estimator of thermal limits when the impedance ($Z = r+ \bm i x$ p.u.) and nominal voltage ($\dot{v}$) is known,
\begin{align}
s^u = \dot{v} \bm e^{-5.0886} \left( \frac{x}{r} \right)^{0.4772}
\end{align}
\noindent
This model is referred to as {\bf TL-Stat}.  The intuition of the model is that the ratio of resistance to impedance provides some insight into the branch's conductor type and configuration, since this ratio should be independent of the line length.

\subsubsection{A Reasonable Upper Bound}

Although the statistical model for branch thermal limits is quite useful, it cannot be applied in cases where data for $r$, $x$, or $\dot{v}$ is missing.  Notable examples include: transformers, where the nominal voltage value differs on both sides of the line, and ideal lines, which do not have an $r$ value.  For these cases, it is helpful to have an alternate method for producing reasonable thermal limits. 

Given that Model \ref{model:acopf} includes reasonable bounds on $|V|$ and $\theta^\Delta$, a reasonable thermal limit can be computed as follows.  For a branch $(l,i,j) \in E$, let $\bm {\theta^{\Delta m}}_{l} = \max(|\bm {\theta^{\Delta l}}_{l}|, |\bm {\theta^{\Delta u}}_{l}|)$. A reasonable value for $\bm {s^u}_{l}$ is
 \begin{align}
(\bm {s^u}_{l})^2 \! = (\bm {v^u}_i)^2 |\bm Y_{l}|^2  ( ({\bm {v^u}_i})^2 \!+ ({\bm {v^u}_j})^2 \!- 2\bm {v^u}_i \bm {v^u}_j\cos(\bm {\theta^{\Delta m}}_{l}) ) \label{eq:thermal:limit}
\end{align}
This model is referred to as {\bf TL-UB} and a more detailed discussion can be found in \cite{nesta}.

%\subsubsection{A Complete Thermal Limit Model}
%
%A robust thermal limit model (TL-Stat) is developed by combining the statistical model and the upper bound model in the following way.  
%A line's thermal limit is computed with both the statistical model and thermal limit upper bound model (TL-UB) and the minimum of these two values
%is used.  If the data required for the statistical model is unavailable on a given line, only thermal limit upper bound is computed.  For test cases that include some line thermal limits, the model is applied only when it produces a tighter limit than the specified value.

%When all of the necessary data is available, the statistical model is used.  In all other cases, the thermal limit upper bound (TL-UB) is used. 

\subsection{Voltage Angle Difference Bounds}

Voltage angle difference bounds are often used as a proxy for capturing voltage angle stability limits on long transmissions lines \cite{9780070359581}.  Unfortunately, no available network datasets include detailed information for these bounds, which are specified in Model \ref{model:acopf}.  To fill this gap all of the \pglibopf~models are given generous voltage angle difference bound of $\bm {\theta^{\Delta l}} \!\!=\!\! -30^\circ$, $\bm {\theta^{\Delta u}} \!\!=\!\! 30^\circ$, which are easily justified by practical voltage stability requirements \cite{9780070359581} and do not impact the best-known solution of any available network.  It is important to emphasize that a voltage angle difference bound of $30^\circ$ is generous enough to be subsumed by the thermal limits provided with all of the networks considered here.  Still, even this generous value has significant implications for the development of power system optimization methods (e.g. \cite{6345338,LPAC_ijoc,qc_opf_tps,Hijazi2017,7763860}).

\section{The \pglibopf~Networks}
\label{sec:pglibopf}

Leveraging the proposed models for completing the missing data in Table \ref{tbl:opf:inst}, the \pglibopf~networks are developed.  Table \ref{tbl:pglib:inst} summarizes which models are used to convert the base network data into \pglibopf~networks.

%Furthermore, it was observed that this voltage angle difference bound and the thermal limit models did not change the optimal solutions of these networks produced by \matpower.
%
%Add note on line loss objective in PEGASE, RTE, SDET, and GOC cases?
%
%Most of the time the \pglibopf~test cases lead to feasible AC-OPF solutions, however, in a few cases it is necessary to override the statistical models to produce AC-OPF feasibility. The special overrides column in Table \ref{tbl:pglib:inst} indicates which networks require additional modifications.  These modifications are detailed as follows:

%\clearpage
\begin{table*}[!ph]
\center
\vspace{-0.5cm}
\caption{\pglibopf~Instance Generation Details}
\footnotesize
%\small
\begin{tabular}{|r||r||c|c|c|c|c|c|r|r|r|r|r|r|r|c|c|}
\hline
 &    & \multicolumn{4}{c|}{Model} &   \\
\pglib & Original & Active &  Reactive   & Gen.  & Thermal  &  Voltage Angle \\
Name &   Name & Gen. & Gen.  & Cost  & Limit & Diff. Bound \\
\hline
%%%%%
\multicolumn{7}{|c|}{Publication Test Cases} \\
\hline
pglib\_opf\_case3\_lmbd & 3-Bus & --- & --- & --- & --- & $30^\circ$ \\
\hline
pglib\_opf\_case5\_pjm & case5 & ---  & --- & --- & TL-Stat & $30^\circ$ \\
\hline
pglib\_opf\_case30\_as & 30 Bus-as & --- & --- & --- & --- & $30^\circ$ \\
\hline
%pglib\_opf\_case30\_fsr & 30 Bus-fsr & --- & --- & --- & --- & $30^\circ$ \\
%\hline
pglib\_opf\_case39\_epri & case39 & --- & --- & --- & --- & $30^\circ$ \\
\hline
%%%%%
\hline
\multicolumn{7}{|c|}{IEEE Power Flow Test Cases}\\
\hline
pglib\_opf\_case14\_ieee &14 Bus & AG-Stat & RG-AM50 & AC-Stat & TL-Stat & $30^\circ$ \\
\hline
pglib\_opf\_case30\_ieee &30 Bus & AG-Stat & RG-AM50 & AC-Stat & TL-Stat & $30^\circ$ \\
\hline
pglib\_opf\_case57\_ieee & 57 Bus & AG-Stat & RG-AM50 & AC-Stat & TL-Stat & $30^\circ$ \\
\hline
pglib\_opf\_case118\_ieee & 118 Bus & AG-Stat & RG-AM50 & AC-Stat & TL-Stat & $30^\circ$ \\
\hline
pglib\_opf\_case300\_ieee & 300 Bus & AG-Stat & RG-AM50 & AC-Stat & TL-UB & $30^\circ$ \\
\hline
\hline
%%%%%
\multicolumn{7}{|c|}{IEEE Dynamic Test Cases}\\
\hline
pglib\_opf\_case162\_ieee\_dtc & 17 Generator & AG-Stat & RG-AM50 & AC-Stat & TL-Stat & $30^\circ$ \\
\hline
\hline
%%%%%
\multicolumn{7}{|c|}{IEEE Reliability Test Systems (RTS)}\\
\hline
pglib\_opf\_case24\_ieee\_rts &  RTS-79 & --- & --- & --- & --- & $30^\circ$ \\
\hline
pglib\_opf\_case73\_ieee\_rts & RTS-96 &--- & --- & --- & --- & $30^\circ$ \\
\hline
\hline
%%%%%
\multicolumn{7}{|c|}{Polish Test Cases}\\
\hline
pglib\_opf\_case2383wp\_k & case2383wp & --- & --- & --- & --- & $30^\circ$ \\
\hline
pglib\_opf\_case2736sp\_k & case2736sp & --- & --- & --- & --- & $30^\circ$ \\
\hline
pglib\_opf\_case2737sop\_k & case2737sop & --- & --- & --- & --- & $30^\circ$ \\
\hline
pglib\_opf\_case2746wop\_k & case2746wop & --- & --- & --- & --- & $30^\circ$ \\
\hline
pglib\_opf\_case2746wp\_k & case2746wp & --- & --- & --- & --- & $30^\circ$ \\
\hline
pglib\_opf\_case3012wp\_k & case3012wp & --- & --- & --- & TL-Stat & $30^\circ$ \\
\hline
pglib\_opf\_case3120sp\_k & case3120sp & --- & --- & --- & TL-Stat & $30^\circ$ \\
\hline
pglib\_opf\_case3375wp\_k & case3375wp & --- & RG-AM50 & --- & TL-Stat & $30^\circ$ \\
\hline
%%%%%
\hline
\multicolumn{7}{|c|}{PEGASE Test Cases}\\
\hline
pglib\_opf\_case89\_pegase & case89pegase & --- & --- & AC-Stat & TL-UB & $30^\circ$ \\
\hline
pglib\_opf\_case1354\_pegase & case1354pegase & --- & --- & AC-Stat & TL-UB & $30^\circ$ \\
\hline
pglib\_opf\_case2869\_pegase & case2869pegase & --- & --- & AC-Stat & TL-UB & $30^\circ$ \\
\hline
pglib\_opf\_case9241\_pegase & case9241pegase & --- & --- & AC-Stat & TL-UB & $30^\circ$ \\
\hline
pglib\_opf\_case13659\_pegase & case13659pegase & --- & --- & AC-Stat & TL-UB & $30^\circ$ \\
\hline
%%%%%
\hline
\multicolumn{7}{|c|}{RTE Test Cases}\\
\hline
pglib\_opf\_case1888\_rte & case1888rte & --- & --- & AC-Stat & TL-UB & $30^\circ$ \\
\hline
pglib\_opf\_case1951\_rte & case1951rte & --- & --- & AC-Stat & TL-UB & $30^\circ$ \\
\hline
pglib\_opf\_case2848\_rte & case2848rte & --- & --- & AC-Stat & TL-UB & $30^\circ$ \\
\hline
pglib\_opf\_case2868\_rte & case2868rte & --- & --- & AC-Stat & TL-UB & $30^\circ$ \\
\hline
pglib\_opf\_case6468\_rte & case6468rte & --- & --- & AC-Stat & TL-UB & $30^\circ$ \\
\hline
pglib\_opf\_case6470\_rte & case6470rte & --- & --- & AC-Stat & TL-UB & $30^\circ$ \\
\hline
pglib\_opf\_case6495\_rte & case6495rte & --- & --- & AC-Stat & TL-UB & $30^\circ$ \\
\hline
pglib\_opf\_case6515\_rte & case6515rte & --- & --- & AC-Stat & TL-UB & $30^\circ$ \\
\hline
%%%%%
\hline
\multicolumn{7}{|c|}{ACTIVSg Test Cases}\\
\hline
pglib\_opf\_case200\_activ & ACTIVSg200 & --- & --- & --- & --- & $30^\circ$ \\
\hline
% pglib\_opf\_case500\_tamu & ACTIVSg500 & --- & --- & --- & --- & $30^\circ$ \\
% \hline
% pglib\_opf\_case2000\_tamu & ACTIVSg2000 & --- & --- & --- & --- & $30^\circ$ \\
% \hline
% pglib\_opf\_case10000\_tamu & ACTIVSg10k & --- & --- & --- & --- & $30^\circ$ \\
% \hline
%%%%%
\hline
\multicolumn{7}{|c|}{Sustainable Data Evolution Technology Test Cases}\\
\hline
pglib\_opf\_case588\_sdet & SDET 500 & --- & --- & AC-Stat & --- & $30^\circ$ \\
\hline
%pglib\_opf\_case2316\_sdet & SDET 2000 & --- & --- & AC-Stat & --- & $30^\circ$ \\
%\hline
pglib\_opf\_case2853\_sdet & SDET 3000 & --- & --- & AC-Stat & --- & $30^\circ$ \\
\hline
pglib\_opf\_case4661\_sdet & SDET 4000 & --- & --- & AC-Stat & --- & $30^\circ$ \\
\hline
%%%%%
\hline
\multicolumn{7}{|c|}{Power Systems Engineering Research Center Test Cases}\\
\hline
pglib\_opf\_case240\_pserc & WECC 240 Bus & --- & --- & AC-Stat & TL-UB & $30^\circ$ \\
\hline
%%%%%
\hline
\multicolumn{7}{|c|}{Grid Optimization Competition Test Cases}\\
\hline
pglib\_opf\_case179\_goc & 179 Bus & --- & --- & AC-Stat & TL-UB & $30^\circ$ \\
\hline
pglib\_opf\_case500\_goc & 02O-173, scenario 72 & --- & --- & --- & --- & $30^\circ$ \\
\hline
pglib\_opf\_case793\_goc & 03O-200, scenario 72 & --- & --- & --- & --- & $30^\circ$ \\
\hline
pglib\_opf\_case2000\_goc & 06O-124, scenario 95 & --- & --- & --- & --- & $30^\circ$ \\
\hline
pglib\_opf\_case2312\_goc & 70O-422, scenario 335 & --- & --- & --- & --- & $30^\circ$ \\
\hline
pglib\_opf\_case2742\_goc & 75O-040, scenario 7 & --- & --- & --- & --- & $30^\circ$ \\
\hline
pglib\_opf\_case3022\_goc & 08O-373, scenario 261 & --- & --- & --- & --- & $30^\circ$ \\
\hline
pglib\_opf\_case3970\_goc & 82O-040, scenario 34 & --- & --- & --- & --- & $30^\circ$ \\
\hline
pglib\_opf\_case4020\_goc & 83O-130, scenario 47 & --- & --- & --- & --- & $30^\circ$ \\
\hline
pglib\_opf\_case4601\_goc & 82O-010, scenario 10 & --- & --- & --- & --- & $30^\circ$ \\
\hline
pglib\_opf\_case4619\_goc & 86O-040, scenario 17  & --- & --- & --- & --- & $30^\circ$ \\
\hline
pglib\_opf\_case4837\_goc & 88O-030, scenario 6 & --- & --- & --- & --- & $30^\circ$ \\
\hline
pglib\_opf\_case4917\_goc & 09O-195, scenario 40 & --- & --- & --- & --- & $30^\circ$ \\
\hline
pglib\_opf\_case9591\_goc & 12O-050, scenario 33 & --- & --- & --- & --- & $30^\circ$ \\
\hline
pglib\_opf\_case10000\_goc & 13O-103, scenario 81 & --- & --- & --- & --- & $30^\circ$ \\
\hline
pglib\_opf\_case10480\_goc & 14O-050, scenario 9 & --- & --- & --- & --- & $30^\circ$ \\
\hline
pglib\_opf\_case19402\_goc & 20O-100, scenario 42 & --- & --- & --- & --- & $30^\circ$ \\
\hline
pglib\_opf\_case24464\_goc & 25O-060, scenario 37 & --- & --- & --- & --- & $30^\circ$ \\
\hline
pglib\_opf\_case30000\_goc & 30O-048, scenario 36 & --- & --- & --- & --- & $30^\circ$ \\
\hline
%%%%%
\end{tabular}
\label{tbl:pglib:inst}
\end{table*}

\subsection{Results}
\label{sec:results}

%
% Tables generated with: Julia v1.0.2, JuMP v0.18.5, PowerModels v0.9.3, Ipopt v0.5.1 and HSL ma27
% Repo: pglib-experiments
% Result: pglib-tbl.py
%

To verify the usefulness of the proposed \pglibopf~networks for benchmarking AC-OPF algorithms, the study from Section \ref{sec:motivation} is revisited.  PowerModels.jl v0.17 \cite{8442948} was used to formulate Model \ref{model:acopf} and its Second-Order Cone (SOC) relaxation \cite{1664986,qc_opf_tps} and both models were solved with IPOPT 3.12 \cite{Ipopt} using the HSL \cite{hsl_lib} linear algebra library on a server with two 2.10GHz Intel CPU and 128GB of RAM.

Table \ref{tbl:gaps_time} presents the results of the base \pglibopf~networks, which are  called the Typical Operating Conditions (TYP) cases.  Many of the optimality gaps have remained small (i.e. below 1\%).  However, a number of networks do exhibit significant optimality gaps, such as pglib\_opf\_case5\_pjm, pglib\_opf\_case30\_ieee,  pglib\_opf\_case162\_ieee\_dtc, pglib\_opf\_case6495\_rte, and pglib\_opf\_case6515\_rte.  This suggests that at least a subset of these cases are interesting for benchmarking AC-OPF algorithms.  

%It may also be worth noting that all of the cases in this table can be solved in less than 100 seconds.

\begin{table*}[!ph]
\small
\center
\caption{AC-OPF Bounds on \pglibopf~TYP Networks.}
\begin{tabular}{|r|r|r||r||r||r|r|r|r|r|r|r|r|r|r|r|r|}
\hline
& & & \$/h & Gap (\%) & \multicolumn{2}{c|}{Runtime (sec.)} \\
Test Case & $|N|$ & $|E|$ & AC & SOC & AC & SOC \\
\hline
\hline
\multicolumn{7}{|c|}{Typical Operating Conditions (TYP)} \\
\hline
pglib\_opf\_case3\_lmbd & 3 & 3 & 5.8126e+03 & 1.32 & $<$1 & $<$1 \\
\hline
pglib\_opf\_case5\_pjm & 5 & 6 & 1.7552e+04 & 14.55 & $<$1 & $<$1 \\
\hline
pglib\_opf\_case14\_ieee & 14 & 20 & 2.1781e+03 & 0.11 & $<$1 & $<$1 \\
\hline
pglib\_opf\_case24\_ieee\_rts & 24 & 38 & 6.3352e+04 & 0.02 & $<$1 & $<$1 \\
\hline
pglib\_opf\_case30\_as & 30 & 41 & 8.0313e+02 & 0.06 & $<$1 & $<$1 \\
\hline
pglib\_opf\_case30\_ieee & 30 & 41 & 8.2085e+03 & 18.84 & $<$1 & $<$1 \\
\hline
pglib\_opf\_case39\_epri & 39 & 46 & 1.3842e+05 & 0.56 & $<$1 & $<$1 \\
\hline
pglib\_opf\_case57\_ieee & 57 & 80 & 3.7589e+04 & 0.16 & $<$1 & $<$1 \\
\hline
pglib\_opf\_case73\_ieee\_rts & 73 & 120 & 1.8976e+05 & 0.04 & $<$1 & $<$1 \\
\hline
pglib\_opf\_case89\_pegase & 89 & 210 & 1.0729e+05 & 0.75 & $<$1 & $<$1 \\
\hline
pglib\_opf\_case118\_ieee & 118 & 186 & 9.7214e+04 & 0.91 & $<$1 & $<$1 \\
\hline
pglib\_opf\_case162\_ieee\_dtc & 162 & 284 & 1.0808e+05 & 5.95 & $<$1 & $<$1 \\
\hline
pglib\_opf\_case179\_goc & 179 & 263 & 7.5427e+05 & 0.16 & $<$1 & $<$1 \\
\hline
pglib\_opf\_case200\_activ & 200 & 245 & 2.7558e+04 & 0.01 & $<$1 & $<$1 \\
\hline
pglib\_opf\_case240\_pserc & 240 & 448 & 3.3297e+06 & 2.78 & 3 & 2 \\
\hline
pglib\_opf\_case300\_ieee & 300 & 411 & 5.6522e+05 & 2.63 & $<$1 & $<$1 \\
\hline
pglib\_opf\_case500\_goc & 500 & 733 & 4.5495e+05 & 0.25 & $<$1 & $<$1 \\
\hline
pglib\_opf\_case588\_sdet & 588 & 686 & 3.1314e+05 & 2.14 & $<$1 & $<$1 \\
\hline
pglib\_opf\_case793\_goc & 793 & 913 & 2.6020e+05 & 1.33 & 2 & $<$1 \\
\hline
pglib\_opf\_case1354\_pegase & 1354 & 1991 & 1.2588e+06 & 1.57 & 4 & 3 \\
\hline
pglib\_opf\_case1888\_rte & 1888 & 2531 & 1.4025e+06 & 2.05 & 8 & 31 \\
\hline
pglib\_opf\_case1951\_rte & 1951 & 2596 & 2.0856e+06 & 0.14 & 18 & 5 \\
\hline
pglib\_opf\_case2000\_goc & 2000 & 3639 & 9.7343e+05 & 0.31 & 7 & 4 \\
\hline
pglib\_opf\_case2312\_goc & 2312 & 3013 & 4.4133e+05 & 1.90 & 6 & 3 \\
\hline
pglib\_opf\_case2383wp\_k & 2383 & 2896 & 1.8682e+06 & 1.04 & 7 & 5 \\
\hline
pglib\_opf\_case2736sp\_k & 2736 & 3504 & 1.3080e+06 & 0.31 & 6 & 4 \\
\hline
pglib\_opf\_case2737sop\_k & 2737 & 3506 & 7.7773e+05 & 0.27 & 5 & 3 \\
\hline
pglib\_opf\_case2742\_goc & 2742 & 4673 & 2.7571e+05 & 1.35 & 28 & 6 \\
\hline
pglib\_opf\_case2746wop\_k & 2746 & 3514 & 1.2083e+06 & 0.37 & 5 & 4 \\
\hline
pglib\_opf\_case2746wp\_k & 2746 & 3514 & 1.6317e+06 & 0.33 & 6 & 5 \\
\hline
pglib\_opf\_case2848\_rte & 2848 & 3776 & 1.2866e+06 & 0.13 & 16 & 7 \\
\hline
pglib\_opf\_case2853\_sdet & 2853 & 3921 & 2.0524e+06 & 0.91 & 9 & 6 \\
\hline
pglib\_opf\_case2868\_rte & 2868 & 3808 & 2.0096e+06 & 0.10 & 14 & 8 \\
\hline
pglib\_opf\_case2869\_pegase & 2869 & 4582 & 2.4628e+06 & 1.01 & 12 & 7 \\
\hline
pglib\_opf\_case3012wp\_k & 3012 & 3572 & 2.6008e+06 & 1.03 & 9 & 12 \\
\hline
pglib\_opf\_case3022\_goc & 3022 & 4135 & 6.0138e+05 & 2.77 & 10 & 5 \\
\hline
pglib\_opf\_case3120sp\_k & 3120 & 3693 & 2.1480e+06 & 0.56 & 9 & 5 \\
\hline
pglib\_opf\_case3375wp\_k & 3374 & 4161 & 7.4382e+06 & 0.55 & 11 & 6 \\
\hline
pglib\_opf\_case3970\_goc & 3970 & 6641 & 9.6099e+05 & 0.34 & 18 & 12 \\
\hline
pglib\_opf\_case4020\_goc & 4020 & 6988 & 8.2225e+05 & 1.23 & 21 & 13 \\
\hline
pglib\_opf\_case4601\_goc & 4601 & 7199 & 8.2624e+05 & 0.54 & 22 & 14 \\
\hline
pglib\_opf\_case4619\_goc & 4619 & 8150 & 4.7670e+05 & 0.91 & 19 & 23 \\
\hline
pglib\_opf\_case4661\_sdet & 4661 & 5997 & 2.2513e+06 & 1.99 & 15 & 11 \\
\hline
pglib\_opf\_case4837\_goc & 4837 & 7765 & 8.7226e+05 & 0.47 & 20 & 11 \\
\hline
pglib\_opf\_case4917\_goc & 4917 & 6726 & 1.3878e+06 & 2.50 & 17 & 9 \\
\hline
pglib\_opf\_case6468\_rte & 6468 & 9000 & 2.0697e+06 & 1.13 & 63 & 27 \\
\hline
pglib\_opf\_case6470\_rte & 6470 & 9005 & 2.2376e+06 & 1.76 & 33 & 23 \\
\hline
pglib\_opf\_case6495\_rte & 6495 & 9019 & 3.0678e+06 & 15.11 & 67 & 24 \\
\hline
pglib\_opf\_case6515\_rte & 6515 & 9037 & 2.8255e+06 & 6.40 & 54 & 22 \\
\hline
pglib\_opf\_case9241\_pegase & 9241 & 16049 & 6.2431e+06 & 2.54 & 48 & 36 \\
\hline
pglib\_opf\_case9591\_goc & 9591 & 15915 & 1.0617e+06 & 0.62 & 60 & 43 \\
\hline
pglib\_opf\_case10000\_goc & 10000 & 13193 & 1.3540e+06 & 1.54 & 47 & 26 \\
\hline
pglib\_opf\_case10480\_goc & 10480 & 18559 & 2.3146e+06 & 1.23 & 65 & 42 \\
\hline
pglib\_opf\_case13659\_pegase & 13659 & 20467 & 8.9480e+06 & 1.39 & 54 & 61 \\
\hline
pglib\_opf\_case19402\_goc & 19402 & 34704 & 1.9778e+06 & 1.19 & 139 & 110 \\
\hline
pglib\_opf\_case24464\_goc & 24464 & 37816 & 2.6295e+06 & 0.97 & 118 & 98 \\
\hline
pglib\_opf\_case30000\_goc & 30000 & 35393 & 1.1423e+06 & 2.89 & 258 & 119 \\
\hline
\end{tabular}
\label{tbl:gaps_time}
\end{table*}

\subsection{Building More Challenging Test Cases}
\label{sec:congested}

The typical operating conditions networks presented in Table \ref{tbl:gaps_time} provide a suitable start for benchmarking AC-OPF algorithms.  However, it is worthwhile to explore if even more challenging test cases can be devised.  To that end, \pglibopf~includes two variants of the base \pglibopf~networks that exhibit even more extreme optimality gaps, the Active Power Increase (API) and Small Angle Difference (SAD) cases.

\paragraph*{Active Power Increase (API) Cases}
It was observed in \cite{pscc_ots, 6345676} that power flow congestion is a key feature of interesting AC Optimal Transmission Switching test cases.  Inspired by this observation, the following Active Power Increase (API) \pglibopf~networks are proposed.  For each of the standard \pglibopf~networks, an optimization problem is solved, which increases the active power demands proportionally throughout the network until the branch thermal limits are binding.  Once a maximal increase in active power demand is determined, the statistical models are applied to appropriately update the other network parameters (e.g. generator capabilities and cost functions).  The results of the API networks are presented in Table \ref{tbl:gaps_time_api}.  As expected, the optimality gaps in these networks have increased significantly with 60\% of cases having optimality gaps above 1\% and eight cases with optimality gaps above 10\%.  This suggests that many of these cases will be useful for benchmarking AC-OPF algorithms.

\begin{table*}[!ph]
\small
\center
\caption{AC-OPF Bounds on \pglibopf~API Networks.}
\begin{tabular}{|r|r|r||r||r||r|r|r|r|r|r|r|r|r|r|r|r|}
\hline
& & & \$/h & Gap (\%) & \multicolumn{2}{c|}{Runtime (sec.)} \\
Test Case & $|N|$ & $|E|$ & AC & SOC & AC & SOC \\
\hline
\hline
\multicolumn{7}{|c|}{Congested Operating Conditions (API)} \\
\hline
pglib\_opf\_case3\_lmbd\_\_api & 3 & 3 & 1.1236e+04 & 9.27 & $<$1 & $<$1 \\
\hline
pglib\_opf\_case5\_pjm\_\_api & 5 & 6 & 7.6377e+04 & 4.09 & $<$1 & $<$1 \\
\hline
pglib\_opf\_case14\_ieee\_\_api & 14 & 20 & 5.9994e+03 & 5.13 & $<$1 & $<$1 \\
\hline
pglib\_opf\_case24\_ieee\_rts\_\_api & 24 & 38 & 1.3494e+05 & 17.88 & $<$1 & $<$1 \\
\hline
pglib\_opf\_case30\_as\_\_api & 30 & 41 & 4.9962e+03 & 44.61 & $<$1 & $<$1 \\
\hline
pglib\_opf\_case30\_ieee\_\_api & 30 & 41 & 1.8044e+04 & 5.46 & $<$1 & $<$1 \\
\hline
pglib\_opf\_case39\_epri\_\_api & 39 & 46 & 2.4967e+05 & 1.72 & $<$1 & $<$1 \\
\hline
pglib\_opf\_case57\_ieee\_\_api & 57 & 80 & 4.9290e+04 & 0.08 & $<$1 & $<$1 \\
\hline
pglib\_opf\_case73\_ieee\_rts\_\_api & 73 & 120 & 4.2263e+05 & 12.87 & $<$1 & $<$1 \\
\hline
pglib\_opf\_case89\_pegase\_\_api & 89 & 210 & 1.3017e+05 & 23.11 & $<$1 & $<$1 \\
\hline
pglib\_opf\_case118\_ieee\_\_api & 118 & 186 & 2.4224e+05 & 29.97 & $<$1 & $<$1 \\
\hline
pglib\_opf\_case162\_ieee\_dtc\_\_api & 162 & 284 & 1.2099e+05 & 4.36 & $<$1 & $<$1 \\
\hline
pglib\_opf\_case179\_goc\_\_api & 179 & 263 & 1.9320e+06 & 9.88 & $<$1 & $<$1 \\
\hline
pglib\_opf\_case200\_activ\_\_api & 200 & 245 & 3.5701e+04 & 0.03 & $<$1 & $<$1 \\
\hline
pglib\_opf\_case240\_pserc\_\_api & 240 & 448 & 4.6406e+06 & 0.67 & 4 & 2 \\
\hline
pglib\_opf\_case300\_ieee\_\_api & 300 & 411 & 6.8499e+05 & 0.85 & $<$1 & $<$1 \\
\hline
pglib\_opf\_case500\_goc\_\_api & 500 & 733 & 6.9241e+05 & 3.44 & 2 & $<$1 \\
\hline
pglib\_opf\_case588\_sdet\_\_api & 588 & 686 & 3.9476e+05 & 1.61 & $<$1 & $<$1 \\
\hline
pglib\_opf\_case793\_goc\_\_api & 793 & 913 & 3.1885e+05 & 13.43 & 2 & $<$1 \\
\hline
pglib\_opf\_case1354\_pegase\_\_api & 1354 & 1991 & 1.4983e+06 & 0.55 & 4 & 3 \\
\hline
pglib\_opf\_case1888\_rte\_\_api & 1888 & 2531 & 1.9539e+06 & 0.23 & 7 & 6 \\
\hline
pglib\_opf\_case1951\_rte\_\_api & 1951 & 2596 & 2.4108e+06 & 0.57 & 8 & 5 \\
\hline
pglib\_opf\_case2000\_goc\_\_api & 2000 & 3639 & 1.4686e+06 & 2.07 & 8 & 4 \\
\hline
pglib\_opf\_case2312\_goc\_\_api & 2312 & 3013 & 5.7152e+05 & 13.09 & 11 & 3 \\
\hline
pglib\_opf\_case2383wp\_k\_\_api & 2383 & 2896 & 2.7913e+05 & 0.01 & 2 & $<$1 \\
\hline
pglib\_opf\_case2736sp\_k\_\_api & 2736 & 3504 & 6.5394e+05 & 10.84 & 7 & 3 \\
\hline
pglib\_opf\_case2737sop\_k\_\_api & 2737 & 3506 & 3.6715e+05 & 5.88 & 6 & 2 \\
\hline
pglib\_opf\_case2742\_goc\_\_api & 2742 & 4673 & 6.4219e+05 & 24.45 & 22 & 5 \\
\hline
pglib\_opf\_case2746wop\_k\_\_api & 2746 & 3514 & 5.1166e+05 & 0.01 & 2 & 2 \\
\hline
pglib\_opf\_case2746wp\_k\_\_api & 2746 & 3514 & 5.8183e+05 & 0.00 & 4 & 2 \\
\hline
pglib\_opf\_case2848\_rte\_\_api & 2848 & 3776 & 1.4970e+06 & 0.25 & 26 & 7 \\
\hline
pglib\_opf\_case2853\_sdet\_\_api & 2853 & 3921 & 2.4578e+06 & 1.96 & 10 & 6 \\
\hline
pglib\_opf\_case2868\_rte\_\_api & 2868 & 3808 & 2.2946e+06 & 0.18 & 23 & 6 \\
\hline
pglib\_opf\_case2869\_pegase\_\_api & 2869 & 4582 & 2.9296e+06 & 1.00 & 13 & 8 \\
\hline
pglib\_opf\_case3012wp\_k\_\_api & 3012 & 3572 & 7.2887e+05 & 0.00 & 4 & 2 \\
\hline
pglib\_opf\_case3022\_goc\_\_api & 3022 & 4135 & 6.5189e+05 & 9.82 & 9 & 4 \\
\hline
pglib\_opf\_case3120sp\_k\_\_api & 3120 & 3693 & 9.3692e+05 & 23.96 & 12 & 4 \\
\hline
pglib\_opf\_case3375wp\_k\_\_api & 3374 & 4161 & 5.8478e+06 & -- & 11 & 352 \\
\hline
pglib\_opf\_case3970\_goc\_\_api & 3970 & 6641 & 1.4557e+06 & 34.58 & 37 & 9 \\
\hline
pglib\_opf\_case4020\_goc\_\_api & 4020 & 6988 & 1.2979e+06 & 17.58 & 39 & 9 \\
\hline
pglib\_opf\_case4601\_goc\_\_api & 4601 & 7199 & 7.9253e+05 & 20.51 & 29 & 10 \\
\hline
pglib\_opf\_case4619\_goc\_\_api & 4619 & 8150 & 1.0299e+06 & 8.30 & 28 & 14 \\
\hline
pglib\_opf\_case4661\_sdet\_\_api & 4661 & 5997 & 2.6953e+06 & 2.64 & 16 & 105 \\
\hline
pglib\_opf\_case4837\_goc\_\_api & 4837 & 7765 & 1.1578e+06 & 8.15 & 36 & 12 \\
\hline
pglib\_opf\_case4917\_goc\_\_api & 4917 & 6726 & 1.5479e+06 & 8.54 & 16 & 8 \\
\hline
pglib\_opf\_case6468\_rte\_\_api & 6468 & 9000 & 2.3135e+06 & 0.82 & 81 & 178 \\
\hline
pglib\_opf\_case6470\_rte\_\_api & 6470 & 9005 & 2.6065e+06 & 1.20 & 52 & 80 \\
\hline
pglib\_opf\_case6495\_rte\_\_api & 6495 & 9019 & 2.9750e+06 & 2.04 & 54 & 22 \\
\hline
pglib\_opf\_case6515\_rte\_\_api & 6515 & 9037 & 3.0617e+06 & 1.89 & 59 & 21 \\
\hline
pglib\_opf\_case9241\_pegase\_\_api & 9241 & 16049 & 7.0112e+06 & 2.67 & 64 & 806 \\
\hline
pglib\_opf\_case9591\_goc\_\_api & 9591 & 15915 & 1.4259e+06 & 13.84 & 145 & 26 \\
\hline
pglib\_opf\_case10000\_goc\_\_api & 10000 & 13193 & 2.3728e+06 & 7.65 & 52 & 27 \\
\hline
pglib\_opf\_case10480\_goc\_\_api & 10480 & 18559 & 2.7627e+06 & 4.78 & 77 & 35 \\
\hline
pglib\_opf\_case13659\_pegase\_\_api & 13659 & 20467 & 9.2842e+06 & 1.85 & 55 & 76 \\
\hline
pglib\_opf\_case19402\_goc\_\_api & 19402 & 34704 & 2.3987e+06 & 4.75 & 195 & 100 \\
\hline
pglib\_opf\_case24464\_goc\_\_api & 24464 & 37816 & 2.4723e+06 & 4.00 & 137 & 151 \\
\hline
pglib\_opf\_case30000\_goc\_\_api & 30000 & 35393 & 1.3530e+06 & 24.79 & 439 & 119 \\
\hline
\end{tabular}
\label{tbl:gaps_time_api}
\end{table*}

\paragraph*{Small Angle Difference (SAD) Cases}
A second approach to modifying the \pglibopf~networks is inspired by recent lines of research \cite{6345338,LPAC_ijoc,qc_opf_tps,Hijazi2017} that indicate voltage angle difference bounds can have significant impacts on power system optimization approaches.  To emphasize these impacts, the following Small Angle Difference (SAD) \pglibopf~networks are proposed.  For each of the standard \pglibopf~networks, an optimization problem is solved in order to find the minimum value of $\bm {\theta^\Delta}$ that can be applied on all of the branches in the network, while retaining a feasible AC power flow.  Once this small value of $\bm {\theta^\Delta}$ is determined, the original test case is updated with this value, which introduces voltage angle difference congestion on the network branches.  The results of the SAD networks are presented in Table \ref{tbl:gaps_time_sad}.  Interestingly, this entirely different approach to modifying the base networks also leads to significant optimality gaps, with 75\% of the SAD networks having optimality gaps above 1\%.  This suggest that many of these cases will be useful for benchmarking AC-OPF algorithms.

It is important to note that in all three result tables, the significant optimality gaps can be caused by two factors: (1) the heuristic for finding feasible AC-OPF solution fails to find the global optimum (e.g. see \cite{6581918}); (2) the SOC convex relaxation is weak (e.g. see \cite{qc_opf_tps}) and does not provide a tight bound on the quality of the AC-OPF solution.  Both factors present interesting avenues for research on AC-OPF algorithms.

\begin{table*}[!ph]
\small
\center
\caption{AC-OPF Bounds on \pglibopf~SAD Networks.}
\begin{tabular}{|r|r|r||r||r||r|r|r|r|r|r|r|r|r|r|r|r|}
\hline
& & & \$/h & Gap (\%) & \multicolumn{2}{c|}{Runtime (sec.)} \\
Test Case & $|N|$ & $|E|$ & AC & SOC & AC & SOC \\
\hline
\hline
\multicolumn{7}{|c|}{Small Angle Difference Conditions (SAD)} \\
\hline
pglib\_opf\_case3\_lmbd\_\_sad & 3 & 3 & 5.9593e+03 & 3.75 & $<$1 & $<$1 \\
\hline
pglib\_opf\_case5\_pjm\_\_sad & 5 & 6 & 2.6109e+04 & 3.62 & $<$1 & $<$1 \\
\hline
pglib\_opf\_case14\_ieee\_\_sad & 14 & 20 & 2.7768e+03 & 21.53 & $<$1 & $<$1 \\
\hline
pglib\_opf\_case24\_ieee\_rts\_\_sad & 24 & 38 & 7.6918e+04 & 9.55 & $<$1 & $<$1 \\
\hline
pglib\_opf\_case30\_as\_\_sad & 30 & 41 & 8.9735e+02 & 7.88 & $<$1 & $<$1 \\
\hline
pglib\_opf\_case30\_ieee\_\_sad & 30 & 41 & 8.2085e+03 & 9.70 & $<$1 & $<$1 \\
\hline
pglib\_opf\_case39\_epri\_\_sad & 39 & 46 & 1.4834e+05 & 0.67 & $<$1 & $<$1 \\
\hline
pglib\_opf\_case57\_ieee\_\_sad & 57 & 80 & 3.8663e+04 & 0.71 & $<$1 & $<$1 \\
\hline
pglib\_opf\_case73\_ieee\_rts\_\_sad & 73 & 120 & 2.2760e+05 & 6.73 & $<$1 & $<$1 \\
\hline
pglib\_opf\_case89\_pegase\_\_sad & 89 & 210 & 1.0729e+05 & 0.73 & $<$1 & $<$1 \\
\hline
pglib\_opf\_case118\_ieee\_\_sad & 118 & 186 & 1.0516e+05 & 8.17 & $<$1 & $<$1 \\
\hline
pglib\_opf\_case162\_ieee\_dtc\_\_sad & 162 & 284 & 1.0869e+05 & 6.48 & $<$1 & $<$1 \\
\hline
pglib\_opf\_case179\_goc\_\_sad & 179 & 263 & 7.6253e+05 & 1.12 & $<$1 & $<$1 \\
\hline
pglib\_opf\_case200\_activ\_\_sad & 200 & 245 & 2.7558e+04 & 0.01 & $<$1 & $<$1 \\
\hline
pglib\_opf\_case240\_pserc\_\_sad & 240 & 448 & 3.4054e+06 & 4.93 & 4 & 2 \\
\hline
pglib\_opf\_case300\_ieee\_\_sad & 300 & 411 & 5.6570e+05 & 2.61 & $<$1 & $<$1 \\
\hline
pglib\_opf\_case500\_goc\_\_sad & 500 & 733 & 4.8740e+05 & 6.67 & $<$1 & $<$1 \\
\hline
pglib\_opf\_case588\_sdet\_\_sad & 588 & 686 & 3.2936e+05 & 6.67 & $<$1 & $<$1 \\
\hline
pglib\_opf\_case793\_goc\_\_sad & 793 & 913 & 2.8580e+05 & 7.97 & 2 & $<$1 \\
\hline
pglib\_opf\_case1354\_pegase\_\_sad & 1354 & 1991 & 1.2588e+06 & 1.57 & 4 & 3 \\
\hline
pglib\_opf\_case1888\_rte\_\_sad & 1888 & 2531 & 1.4139e+06 & 2.82 & 9 & 19 \\
\hline
pglib\_opf\_case1951\_rte\_\_sad & 1951 & 2596 & 2.0924e+06 & 0.46 & 17 & 5 \\
\hline
pglib\_opf\_case2000\_goc\_\_sad & 2000 & 3639 & 9.9288e+05 & 1.52 & 7 & 4 \\
\hline
pglib\_opf\_case2312\_goc\_\_sad & 2312 & 3013 & 4.6235e+05 & 3.93 & 7 & 4 \\
\hline
pglib\_opf\_case2383wp\_k\_\_sad & 2383 & 2896 & 1.9112e+06 & 2.86 & 8 & 5 \\
\hline
pglib\_opf\_case2736sp\_k\_\_sad & 2736 & 3504 & 1.3266e+06 & 1.58 & 7 & 5 \\
\hline
pglib\_opf\_case2737sop\_k\_\_sad & 2737 & 3506 & 7.9095e+05 & 1.88 & 7 & 4 \\
\hline
pglib\_opf\_case2742\_goc\_\_sad & 2742 & 4673 & 2.7571e+05 & 1.36 & 28 & 5 \\
\hline
pglib\_opf\_case2746wop\_k\_\_sad & 2746 & 3514 & 1.2337e+06 & 2.32 & 7 & 4 \\
\hline
pglib\_opf\_case2746wp\_k\_\_sad & 2746 & 3514 & 1.6669e+06 & 2.19 & 7 & 5 \\
\hline
pglib\_opf\_case2848\_rte\_\_sad & 2848 & 3776 & 1.2890e+06 & 0.25 & 16 & 6 \\
\hline
pglib\_opf\_case2853\_sdet\_\_sad & 2853 & 3921 & 2.0692e+06 & 1.70 & 9 & 6 \\
\hline
pglib\_opf\_case2868\_rte\_\_sad & 2868 & 3808 & 2.0213e+06 & 0.59 & 16 & 6 \\
\hline
pglib\_opf\_case2869\_pegase\_\_sad & 2869 & 4582 & 2.4687e+06 & 1.12 & 11 & 7 \\
\hline
pglib\_opf\_case3012wp\_k\_\_sad & 3012 & 3572 & 2.6195e+06 & 1.56 & 10 & 14 \\
\hline
pglib\_opf\_case3022\_goc\_\_sad & 3022 & 4135 & 6.0143e+05 & 2.77 & 9 & 5 \\
\hline
pglib\_opf\_case3120sp\_k\_\_sad & 3120 & 3693 & 2.1749e+06 & 1.52 & 11 & 6 \\
\hline
pglib\_opf\_case3375wp\_k\_\_sad & 3374 & 4161 & 7.4382e+06 & 0.55 & 11 & 7 \\
\hline
pglib\_opf\_case3970\_goc\_\_sad & 3970 & 6641 & 9.6555e+05 & 0.76 & 20 & 12 \\
\hline
pglib\_opf\_case4020\_goc\_\_sad & 4020 & 6988 & 8.8969e+05 & 8.66 & 28 & 13 \\
\hline
pglib\_opf\_case4601\_goc\_\_sad & 4601 & 7199 & 8.7818e+05 & 6.37 & 29 & 14 \\
\hline
pglib\_opf\_case4619\_goc\_\_sad & 4619 & 8150 & 4.8435e+05 & 2.00 & 22 & 15 \\
\hline
pglib\_opf\_case4661\_sdet\_\_sad & 4661 & 5997 & 2.2610e+06 & 1.96 & 15 & 12 \\
\hline
pglib\_opf\_case4837\_goc\_\_sad & 4837 & 7765 & 8.7712e+05 & 0.85 & 22 & 11 \\
\hline
pglib\_opf\_case4917\_goc\_\_sad & 4917 & 6726 & 1.3890e+06 & 2.52 & 17 & 9 \\
\hline
pglib\_opf\_case6468\_rte\_\_sad & 6468 & 9000 & 2.0697e+06 & 1.12 & 63 & 30 \\
\hline
pglib\_opf\_case6470\_rte\_\_sad & 6470 & 9005 & 2.2416e+06 & 1.91 & 35 & 22 \\
\hline
pglib\_opf\_case6495\_rte\_\_sad & 6495 & 9019 & 3.0678e+06 & 15.11 & 69 & 23 \\
\hline
pglib\_opf\_case6515\_rte\_\_sad & 6515 & 9037 & 2.8698e+06 & 7.85 & 57 & 21 \\
\hline
pglib\_opf\_case9241\_pegase\_\_sad & 9241 & 16049 & 6.3185e+06 & 2.47 & 52 & 36 \\
\hline
pglib\_opf\_case9591\_goc\_\_sad & 9591 & 15915 & 1.1674e+06 & 9.56 & 83 & 42 \\
\hline
pglib\_opf\_case10000\_goc\_\_sad & 10000 & 13193 & 1.4902e+06 & 6.77 & 44 & 26 \\
\hline
pglib\_opf\_case10480\_goc\_\_sad & 10480 & 18559 & 2.3147e+06 & 1.23 & 65 & 41 \\
\hline
pglib\_opf\_case13659\_pegase\_\_sad & 13659 & 20467 & 9.0422e+06 & 1.68 & 62 & 43 \\
\hline
pglib\_opf\_case19402\_goc\_\_sad & 19402 & 34704 & 1.9838e+06 & 1.49 & 162 & 108 \\
\hline
pglib\_opf\_case24464\_goc\_\_sad & 24464 & 37816 & 2.6540e+06 & 1.84 & 145 & 101 \\
\hline
pglib\_opf\_case30000\_goc\_\_sad & 30000 & 35393 & 1.2866e+06 & 11.84 & 336 & 120 \\
\hline
\end{tabular}
\label{tbl:gaps_time_sad}
\end{table*}

\section{Conclusions}
\label{sec:conclusion}

This report has highlighted some of the shortcomings of using existing network datasets for benchmarking AC Optimal Power Flow algorithms and has proposed the \pglibopf~networks to mitigate them.  Leveraging data-driven models, \pglibopf~ensures that all of the networks have reasonable values for key power network parameters, including generation injection limits, generation costs, and branch thermal limits.  Furthermore, the active power increase and small angle difference network variants are developed to provide additional challenging cases for benchmarking.  A detailed validation study demonstrates that the majority of the \pglibopf~networks exhibit significant optimality gaps and are therefore useful for benchmarking AC-OPF algorithms.

It is important to emphasize that while the primary challenge of this work has been to develop realistic and challenging network datasets for benchmarking the AC-OPF problem presented in Model \ref{model:acopf}, there is still a significant gap in the models used for industry-grade AC optimal power flow studies \cite{stott1987, real_opf,capitanescu2016}.  Some of the key extensions include: information about configurable assets such as bus shunts, switches, and transformer tap settings; N-1 contingency cases;  branch thermal limits for short- and long-term overloading; and generator capability curves.  As the research community is able to addresses the challenges presented in the \pglibopf~networks, it is important to also consider more realistic extensions of Model \ref{model:acopf} and to curate new \pglib~repositories for those model variants. 

Although \pglibopf~has highlighted some advantages for AC-OPF benchmarking, its network datasets are still by-in-large synthetically generated.  This highlights the continued need for industry engagement in the development of more detailed and realistic network datasets for benchmarking AC-OPF algorithms.  We hope that the \pglibopf~networks proposed herein will be sufficient for the research community to benchmark AC-OPF algorithms, while more real-world network datasets are developed and contributed to the \pglibopf~repository in the years to come.

\appendix[Extensions and Alternative Applications]

While this report focuses on the AC-OPF formulation described in Model~\ref{model:acopf}, there are a variety of possible modifications, extensions, and other problems that are relevant to many power system researchers. This appendix first presents common modifications of Model~\ref{model:acopf} and then summarizes several possible alternative uses of the \pglib~data.

\subsection{Current Flow Limits}

The OPF formulation in Model~\ref{model:acopf} considers line-flow limits that are based on apparent power~\eqref{eq:acopf:7} and voltage angle differences~\eqref{eq:acopf:8}. A common modified problem formulation limits the magnitudes of the current flows on each line. Specifically, the following constraints either replace or augment the apparent power flow limits in~\eqref{eq:acopf:7}:
\begin{subequations}
\begin{align}
& I_{lij} = \left( \bm Y_{l} + \bm i\frac{\bm {b^c}_{l}}{2} \right) \frac{V_i}{|\bm{T}_{l}|^2} - \bm Y_{l} \frac{V_j}{\bm{T}^*_{l}} \;\; \forall (l,i,j) \in E \\
& I_{lji} = \left( \bm Y_{l} + \bm i\frac{\bm {b^c}_{l}}{2} \right) V_j - \bm Y_{l} \frac{V_i}{\bm{T}_{l}} \;\; \forall (l,i,j) \in E \\
& |I_{lij}| \leq {\bm I}_{l}^{\bm u} \;\; \forall (l,i,j) \in E \cup E^R
\end{align}
\end{subequations}
Note that the per unit value of ${\bm I}_{l}^{\bm u}$ is often assumed to be equivalent to the per unit value of $\bm {s^u}_{l}$.

% S_{lij} = \left( \bm Y^*_{l} - \bm i\frac{\bm {b^c}_{l}}{2} \right) \frac{|V_i|^2}{|\bm{T}_{l}|^2} - \bm Y^*_{l} \frac{V_i V^*_j}{\bm{T}_{l}} \;\; \forall (l,i,j)\in E \label{eq:acopf:5} \\
% & S_{lji} = \left( \bm Y^*_{l} - \bm i\frac{\bm {b^c}_{l}}{2} \right) |V_j|^2 - \bm Y^*_{l} \frac{V^*_i V_j}{\bm{T}^*_{l}}

% Write out current flow expression and explain data format

% \subsection{Multi-Period OPF Formulations}
% The OPF formulation in Model~\ref{model:acopf} considers a single period of the power system behavior. More general ``multi-period'' formulations couple multiple periods where quantities such as the load demands vary over time. Consecutive periods are coupled by generator ramp rates, the dynamics of energy storage devices and shiftable loads, etc.

% citation to multiperiod formulations?

% planning problems with differing load scenarios. Some test cases (e.g., Polish systems) provide this information via different snapshots. Can also use RTS system data to construct a year of hourly data regarding loading scenarios. Citations?

\subsection{Branch Charging Model}

The $\Pi$-circuit branch model of this work (i.e., Figure~\ref{fig:linemodel}), only considers the susceptance impacts of line changing.  After the inception of that \matpower~data format, a number of commercial power flow tools now support the following, more general model, of line charging:
\begin{subequations}
\begin{align}
& S_{lij} = \left( \bm Y_{l} + \bm Y^c_{lij} \right)^* \frac{|V_i|^2}{|\bm{T}_{l}|^2} - \bm Y^*_{l} \frac{V_i V^*_j}{\bm{T}_{l}} \;\; \forall (l,i,j) \in E \\
& S_{lji} = \left( \bm Y_{l} + \bm Y^c_{lji} \right)^* |V_j|^2 - \bm Y^*_{l} \frac{V^*_i V_j}{\bm{T}^*_{l}} \;\; \forall (l,i,j) \in E
\end{align}
\end{subequations}
where the values $\bm Y^c_{lij}$ and $\bm Y^c_{lji}$ represent the line charging admittance on the {\em from} and {\em to} sides of the branch respectively.  This model generalizes the branch model from Figure \ref{fig:linemodel} by incorporating line charge conductance effects and asymmetrical charging effects.

\subsection{Generator Capability Curves}
Model~\ref{model:acopf} represents generators with box constraints that independently limit active and reactive power outputs. A more detailed model recognizes that the current flows inside of a generator are jointly dependent on both the active and reactive power outputs. The generator must be operated to limit the heating caused by these internal current flows. Thus, the more detailed ``capability curve'' generator model (also known as a ``D-curve'' model) forms generator limits that couple the active and reactive power outputs~\cite{jackson1971}. Estimates of the additional data needed for the generator capability curve model can be extracted from the box constraints on active and reactive power injections provided in typical datasets~\cite{molzahn2015naps,park2017curves}. Additionally, the \matpower{} data format is capable of representing a piecewise-linear approximation of the generator capability curve model~\cite{matpower}.

\subsection{Voltage \& Reactive Power Control}
In practice, reactive power management devices, such as bus shunts and transformer taps, play a critical role in managing the voltage profile of an AC power network and improving power quality in congested parts of a network \cite{7286009}.  A notable limitation of the OPF formulation presented in Model~\ref{model:acopf} is the limited amount of reactive power control devices, which may bias the feasible solutions to a specific voltage profile provided with the network data.  At this time, extensions of the Model~\ref{model:acopf} that consider reactive power controls are an active area of research with promising initial results \cite{GOCBench,pes_opf_comp,7038395,h_minlp_opf}. However, a variety of subtle challenges arise when incorporating such devices into an OPF problem formulation \cite{real_opf} and continued research is required to arrive at a broadly accepted generalization of Model~\ref{model:acopf} that incorporates more reactive control devices.

\subsection{Other Problem Formulations}
The OPF formulation in Model~\ref{model:acopf} considers a single period of steady-state power system behavior. Power system engineers solve a wide variety of other optimization and control problems relevant to the design and operation of power systems~\cite{sasson1974,momoh2008,wood2013}. Formulating many of these problems requires augmenting the data available in \pglibopf{} with other information. Table~\ref{tbl:other_problems} summarizes several other classes of optimization and control problems, indicates examples of information that may be required for these problems beyond the data provided in \pglibopf{}, and provides selected references for each class of problems. 
% Also note that many of the problems described in Table~\ref{tbl:other_problems} can be extended to variants that consider, e.g., stochasticity~\cite{bienstock2014,tahanan2015,wang2018} and distributed formulations~\cite{kekatos2017,molzahn_dorfler_sandberg_low_chakrabarti_baldick_lavaei-distributed_optimization_control}. 
While currently beyond the scope of this effort, extension of \pglib{} to incorporate the data necessary for test cases that are applicable to these and other classes of problems is an important topic for future work.

\begin{table*}[tb]
\caption{Other Problem Formulations}
{\renewcommand{\arraystretch}{1.2}%
\footnotesize
\label{tbl:other_problems}
\centering
\begin{tabular}{|M{0.12\textwidth}|M{0.3\textwidth}|M{0.42\textwidth}|M{0.08\textwidth}|}
\hline 
\multicolumn{1}{|c|}{{Problem}} & \multicolumn{1}{c|}{{Description}} & \multicolumn{1}{|c|}{{Additional Information Required}} & \multicolumn{1}{|c|}{{Refs.}} \\ \hline\hline
Power Flow & Determine a voltage profile corresponding to a specified set of power injections and voltage magnitudes. &  
\begin{itemize}[leftmargin=*]
\item Bus type specifications (PV, PQ, Slack) for each bus.
\item Active power injections for PV and PQ buses.
\item Reactive power injections for PQ buses.
\item Voltage magnitudes for the slack bus and PV buses.
\end{itemize} & \cite{stott1974,acc2016_tutorial} \tabularnewline[3.5em]\hline
Multi-Period OPF & Extend the single-period OPF formulation of Model~\ref{model:acopf} to multiple periods with time-varying load demands and energy storage models. & 
\begin{itemize}[leftmargin=*]
\item Time-varying load demands and renewable generation.
\item Generator ramp rates.
\item Models of energy storage devices.
\end{itemize} & \cite{rts96,alguacil2000} \tabularnewline[2.5em]\hline
Unit Commitment & Schedule generator on/off statuses and power outputs considering costs and constraints associated with operating, starting up, and shutting down generators. & 
\begin{itemize}[leftmargin=*]
\item Time-varying load demands and renewable generation.
\item Generator start-up and shut-down costs, ramp rates, minimum run times, and minimum down times.
\item Models of energy storage devices.
\end{itemize} & \cite{sen1998,padhy2004,tahanan2015,knueven2018} 
% \cite{sen1998},\cite{padhy2004},\newline\hspace*{0pt} \cite{tahanan2015},\cite{knueven2018} 
\tabularnewline[2.5em]\hline
Network \\ Reconfiguration & Optimize the topology of a transmission or distribution network by changing the statuses of switches. &  
\begin{itemize}[leftmargin=*]
\item Locations of switches.
\end{itemize} & \cite{4492805,hedman2011},\newline\hspace*{0pt}\cite{hijazi2014reconfiguration} \tabularnewline[1.5em]\hline
State Estimation & Determine the operating point that minimizes a weighted error metric for a given set of measurements. &  
\begin{itemize}[leftmargin=*]
\item Locations and values for the measurements.
\item Characterizations of measurement noise.
\end{itemize} & \cite{wu1990,abur04a,kekatos2017} \tabularnewline[2.5em]\hline
Voltage Stability Analyses & Compute various margins to the system's loadability limit. &  
\begin{itemize}[leftmargin=*]
\item Models for variations in the load demands and generator outputs.
\item Models for the actions of voltage control devices such as switched shunt capacitors and tap-changing transformers.
\end{itemize} & \cite{taylor94a,ieee02a} \tabularnewline[2.5em]\hline
Small-Signal Dynamic Stability Analyses & Simulate or analytically characterize the behavior of a power system subject to small perturbations. &  
\begin{itemize}[leftmargin=*]
\item Dynamical models for devices such as generators and loads, often represented as a system of differential-algebraic equations.
\item Models for the perturbations.
\end{itemize} & \cite{sauer98a} \tabularnewline[2.5em]\hline
Transient Stability Analyses & Simulate or analytically characterize the behavior of a power system following a large disturbance. &  
\begin{itemize}[leftmargin=*]
\item Dynamical models for devices such as generators and loads, often represented as a system of differential-algebraic equations.
\item Models for the disturbances.
\end{itemize} & \cite{sauer98a,abhyankar2017tsopfI,abhyankar2017tsopfII} \tabularnewline[2.5em]\hline
Cascading Failure Analyses & Simulate and analyze the processes by which failures of certain system devices can lead to cascades of subsequent failures. &  
\begin{itemize}[leftmargin=*]
\item Models relating overload amounts to the failure probabilities of each device.
\item System responses to failures of various devices.
\end{itemize} & \cite{vaiman2012} \tabularnewline[2.5em]\hline
Reliability \mbox{Analyses} and Expansion Planning & Compute estimates of reliability metrics such as loss of load probability. Identify investments in new devices such as generators and transmission lines that facilitate reliable operation. &  
\begin{itemize}[leftmargin=*]
\item Multiple snapshots of load demands representing time-series or scenario data.
\item Failure rates for various devices such as generators and lines.
\end{itemize} & \cite{6407493,rts96},\newline\hspace*{0pt}\cite{romero2002} \tabularnewline[2.5em]\hline
Coupled Transmission / Distribution System Analyses & Analyses that jointly consider models of transmission and distribution systems. &  
\begin{itemize}[leftmargin=*]
\item Distribution system models, possibly with three-phase unbalanced network representations.
\item Model for the couplings between the transmission and distribution systems.
\end{itemize} & \cite{aristidou2014,palmintier2016,abhyankar2017} \tabularnewline[3.5em]\hline
Coupled Infrastructure Analyses & Analyses of electric power systems coupled with other networks such as natural gas, water, transportation, and communication systems. &  
\begin{itemize}[leftmargin=*]
\item Models for the other infrastructures and their couplings with the electric grid.
\end{itemize} & \cite{khaitan2013,zlotnik2017,zamzam2018} \tabularnewline[2.5em]\hline
Stochastic \\ Formulations & Problem variants that consider uncertainty, often with respect to the power injections. &  
\begin{itemize}[leftmargin=*]
\item Models of the uncertainty characteristics.
\item Recourse models describing how each device responds to the uncertainty realizations.
\end{itemize} & \cite{wang2018,tahanan2015},\newline\hspace*{0pt}\cite{bienstock2014}\tabularnewline[1.5em]\hline
Distributed \\ Formulations & Problem variants that use distributed rather than centralized formulations. &  
\begin{itemize}[leftmargin=*]
\item Models of the communications among the distributed agents that work together to solve the problem.
\end{itemize} & \cite{kekatos2017,molzahn_dorfler_sandberg_low_chakrabarti_baldick_lavaei-distributed_optimization_control}\tabularnewline[1.5em]\hline
Contingency-Constrained Formulations & Problem variants that ensure the security of system operations after the occurrence of certain contingencies. &  
\begin{itemize}[leftmargin=*]
\item A list of contingencies.
\item Recourse models describing how each device responds to the contingency.
\end{itemize} & \cite{stott1987,real_opf,capitanescu2016}\tabularnewline[1.5em]\hline
\end{tabular}%
}%
\end{table*}

%abhyankar2017tsopfI,abhyankar2017tsopfII

% hedman2011

% Table for power flow, unit commitment, state estimation, transmission switching, contingency analysis, expansion planning, transient stability, stochastic optimization, 
% describe what additional data would be required to formulate these problems

% question about power flow: what are the data provided in the Pg, Qg, Vg fields?

%\clearpage
\bibliographystyle{IEEEtran}
\bibliography{refrences.bib}

LA-UR-18-29054

\end{document}